\documentstyle[12pt,amstex,amssymb]{article}
\newtheorem{thm}{Theorem}[section]
\newtheorem{prop}[thm]{Proposition}

\newtheorem{cor}[thm]{Corollary}

\def\p{\partial}
\def\w{\wedge}
\def\ker{\operatorname{ker}}

\def\dim{\operatorname{dim}}

\def\rk{\operatorname{rk}}

\newcommand{\pf}{\noindent{\em Proof.}\ }

\newcommand{\CC}{{\Bbb C}}

\newcommand{\ZZ}{{\Bbb Z}}

\newcommand{\lan}{\langle}
\newcommand{\ran}{\rangle}

\renewcommand{\a}{\alpha}
\renewcommand{\b}{\beta}
\newcommand{\g}{\gamma}
\renewcommand{\d}{\delta}
\newcommand{\e}{\varepsilon}
\renewcommand{\l}{\lambda}
\newcommand{\m}{\mu}
\renewcommand{\o}{\omega}

\renewcommand{\t}{\tau}
\newcommand{\G}{\Gamma}
\newcommand{\D}{\Delta}
\renewcommand{\L}{\Lambda}

\newcommand{\cald}{{\cal D}}

\newcommand{\calf}{{\cal F}}

\newcommand{\calh}{{\cal H}}

\newcommand{\calm}{{\cal M}}

\newcommand{\calo}{{\cal O}}

\begin{document}

\setlength{\baselineskip}{18pt}

\author{ Alexander B. GIVENTAL
\thanks{Supported by Alfred P. Sloan Foundation
and by NSF Grant DMS--9321915.} }

\title{EQUIVARIANT GROMOV-WITTEN INVARIANTS}

\date{Department of Mathematics\\University of California\\
 Berkeley, California 94720-3840 USA }

\maketitle

The objective of this paper is to describe some construction and applications
of the equivariant counterpart to the Gromov-Witten (GW) theory, i.e.
intersection theory on spaces of (pseudo-) holomorphic curves in (almost-)
Kahler manifolds.

Given a Killing action of a compact Lie group $G$ on a compact Kahler manifold
$X$, the equivariant GW-theory provides, as we will show in Section $3$,
the equivariant cohomology space
$H^*_G(X)$ with a {\em Frobenius structure} (see \cite{Db}).
We discuss applications of the equivariant theory to the computation
(\cite{GK},\cite{K}) of quantum cohomology algebras of flag manifolds
(Section $5$), to the
simultaneous diagonalization of the quantum cup-product operators
(Sections $7$,$8$), and to the
$S^1$-equivariant Floer homology theory on the loop space $LX$ (see Section $6$
and \cite{HG},\cite{HG1}).

In Sections $9$ --- $11$ we combine the general theory developed in Sections
$1$ --- $6$ with the fixed point localization technique \cite{Kn} in order to
prove the mirror conjecture (in the form suggested in \cite{HG}) for
projective complete intersections. We manage to interpret in rigorous terms
of stable maps \cite{Kn} all the
results from \cite{HG1, HG} which were motivated by geometry of loop spaces
of projective complete intersections.
In particular, we prove in Section $11$ the following

\bigskip

{\bf Theorem.} {\em
Consider the Picard-Fuchs differential equation
\[ (\frac{d}{dt})^4 I =5e^t(5\frac{d}{dt}+1)(5\frac{d}{dt}+2)(5\frac{d}{dt}+3)
(5\frac{d}{dt}+4) I \]
satisfied by the periods
\[ I(t) =
\int_{\g_t ^3 } \frac{du_0\w ... \w du_4 }{d(u_0+...+u_4)\w d(u_0...u_4)} \]
of the non-vanishing holomorphic $3$-forms on the Calabi-Yau $3$-folds $Y_t$
with Hodge numbers $h^{2,1}=1,\ h^{1,1}=101$ given by the affine equations
$Y_t: u_0+...+u_4=1, u_0...u_4=e^t$.

Pick the basis $I_0,...,I_3$ of solutions to this differential equation
determined by
\[ I_0(t)+I_1(t)P+I_2(t)P^2+I_3(t)P^3 = \
\sum _{d=0}^{\infty } e^{(P+d)t} \frac{\Pi _{m=1}^{5d}(5P+m)}
{\Pi _{m=1}^d (P+m)^5} \ \ (\ \text{mod}\ P^4) .\]

Introduce the new variable $T(t)=I_1(t)/I_0(t)$.

Then

\[ \frac{I_0}{I_0}+\frac{I_1}{I_0}(t(T))P+
\frac{I_2}{I_0}(t(T))P^2+\frac{I_3}{I_0}(t(T))P^3 \
=  \]
\[ e^{PT} +  \frac{P^2}{5} \sum_{d=1}^{\infty} n_d d^3 \sum_{k=1}^{\infty}
\frac{e^{(P+kd)T}}{(P+kd)^2} \ \ (\ \text{mod}\ P^4) ,\]
where the components of the RHS form the basis of solutions to the differential
equation
\[ (\frac{d}{dT})^2 \frac{1}{K(e^T)} (\frac{d}{dT})^2 \ J \ = 0
\ \ \text{with} \ \ K(q)=5+\sum_{d=1}^{\infty } n_d d^3 \frac{q^d}{1-q^d} ,\]
and $n_d$ is the virtual number of degree $d$ rational curves in $\CC P^4$
situated on a generic degree $5$ hypersurface $X$, a Calabi-Yau $3$-fold with
Hodge numbers $h^{2,1}=101,\ h^{1,1}=1$.

Analogous results hold for any non-singular Calabi-Yau $3$-dimensional
projective complete intersection $X$.
 }

\bigskip

The virtual numbers of rational curves on a Calabi-Yau $3$-fold $X$ are defined
in several equivalent ways in the quantum cohomology theory (see \cite{AM, M})
and are equal to the algebraic number of such curves on $X$ provided with a
generic almost Kahler structure. It is known that for generic quintic
hypersurfaces $X\subset \CC P^4$ the virtual number $n_d$ coincide with
the number of the degree $d$ rational curves in $\CC P^4$ situated in $X$
at least for $d\leq 24$.
The number $n_1=2875$ of straight lines on a generic quintic $3$-fold has
been known since the last century, $n_2=609250$ and $n_3=317206375$ were found
(see \cite{COGP}) several years ago, while $n_4=242467530000$ was predicted in
\cite{COGP} and confirmed in \cite{Kn} (as an illustration of a method that
allows in principle to find each $n_d$).
The simultaneous description of all the numbers $n_d$ given in the theorem
was conjectured in \cite{COGP} on the basis of physical ideas of mirror
symmetry between the Calabi-Yau manifolds $X$ and $Y$ whose Hodge diamonds
happened to be mirror-symmetric to one another.

As far as we know, our Theorem and its generalization to Calabi-Yau
projective complete intersections given in Section $11$ provide the first
examples of Calabi-Yau manifolds for which predictions of the mirror symmetry
are verified for rational curves of all degrees.

The results of Sections $9$ -- $11$ can be immediately carried over to complete
intersections in products of projective spaces. The method can be also applied
to complete intersections in general toric varieties where however some
generalization of our algebraic formalism and some refinement in foundations
of the equivariant Gromov -- Witten theory would be necessary.

\bigskip

I am thankful to S. Barannikov, I. Grojnowski, B. Kim, D. Morrison,
R. Plesser, A. Schwartz, A. Varchenko for numerous stimulating discussions
and especially to M. Kontsevich who taught me his approach to Gromov --
Witten theory. The proof of the theorem formulated above has
grown out of our joint attempt in Spring $95$ to prove it using the method
\cite{Kn} of summation over trees. The influence of our discussions on
other results of this paper is also significant.

\section{Moduli spaces of stable maps}
\label{sec1}
\setcounter{equation}{0}

It was M. Gromov \cite{Gr} who first
suggested to construct (and constructed some) topological invariants of
a symplectic manifold $X$
as bordism classes of spaces of pseudo-holomorphic curves in $X$. Recently
M. Kontsevich \cite{Kn} suggested the concept of {\em stable maps} which gives
rise to an adequate compactification of these spaces. We recall here some
basic facts from \cite{Kn} about these compactifications.

Let $(C,p)$ be a compact connected complex curve with only double singular
points
and with $n$ ordered non-singular {\em marked points} $(p_1,...,p_n)$.
Two holomorphic maps
$(C, p) \to X, \ (C', p') \to X$ to an almost-Kahler manifold $X$ are called
{\em equivalent} if they can be identified by a holomorphic isomorphism
$(C, p) \to (C', p')$. A holomorphic map $(C,p)\to X$ is called {\em stable}
if it does not have infinitesimal automorphisms. In other words, a map is
unstable if either it is constant on a genus $0$ irreducible component of
$C$ with $< 3$ {\em special} ($=$ marked or singular) points or
if $C$ is a torus, carries no marked points and the map is constant.
A stable map may have a non-trivial finite automorphism group.

According to Gromov's compactness theorem \cite{Gr}, any sequence of
holomorphic maps
$C\to X$ of a nonsingular compact curve $C$ has a subsequence
Hausdorff-convergent to a holomorphic map $\hat C\to X$ of (may be reducible)
curve $\hat C$ of the same genus $g$ and representing the same total homology
class $d\in H_2(X,\Bbb Z)$. A refinement of this theorem from \cite{Kn} says
that equivalence classes of stable maps $C\to X$ with given $g,n,d$ form a
single
compact Hausdorff space --- the moduli space of stable maps --- which we
denote $X_{g,n,d}$. Here $g=\dim H^1(C,\calo)=1-
\chi (C \backslash C^{\mbox{sing}})$.

In the case $X=pt$ the moduli spaces coincide with Deligne-Mumford
compactifications $\calm _{g,n}$ of moduli spaces of genus $g$ Riemannian
surfaces with $n$ marked points. They are compact nonsingular orbifolds
(i.e. local
quotients of nonsingular manifolds by finite groups) and thus bear the rational
fundamental cycle which allows one to build up intersection theory.
In general, the moduli spaces $X_{g,n,d}$ are singular and may have ``wrong"
dimension, and the idea of the program started in \cite{KM, Kn} is to provide
$X_{g,n,d}$ with virtual fundamental cycles insensitive to perturbations of the
almost-Kahler structure on $X$.
In some nice cases however the spaces $X_{g,n,d}$ are already nonsingular
orbifolds of the ``right" dimension.

A compact complex manifold is called {\em
ample} if it is a homogeneous space of its Lie algebra of holomorphic
vector fields.

\begin{thm}[\cite{Kn, BM}]
   If $X$ is ample then all non-empty moduli spaces
$X_{n,d}$   of genus $0$   stable maps are compact nonsingular complex
orbifolds of ``right" dimension  $\lan c_1(T_x),d\ran +\dim_{\Bbb C}X+n-3$.
\end{thm}

Additionally, there are canonical morphisms $X_{n,d}\to X_{n-1,d}, \
X_{n,d}\to\calm _{0,n},\  X_{n,d}\to X^n$ between the moduli spaces $X_{n,d}$
called {\em forgetful, contraction} and {\em evaluation} (and defined
by forgetting one of the marked points, forgetting the map and evaluating
the map at marked points respectively). We refer to \cite{Kn, BM} for details
of their construction.

In the rest of this paper we will stick to ample manifolds; we
comment however on  which results are expected to hold in greater generality.
A number of recent preprints by B. Behrend -- B. Fantechi, J.Li -- G. Tian,
T. Fukaya -- K. Ono shows that Kontsevich's ``virtual fundamental cycle''
program is being realized successfully and leaves no doubts that these
generalizations are correct. Still some verifications are necessary in order
make them precise theorems.

\section{Equivariant correlators}
\label{sec2}
\setcounter{equation}{0}

The Gromov-Witten theory borrows from quantum field theory the name {\em
(quantum) correlators} for numerical topological characteristics of the moduli
spaces $X_{n,d}$ (characteristic numbers) and borrows from bordism theory
the construction of such correlators as integrals
of suitable wedge-products of various universal cohomology classes
(characteristic classes of the GW theory) over the fundamental cycle.

We list here some such characteristic classes.

\begin{enumerate}
\item Pull-backs of cohomology classes from $X^n$ by the evaluation
maps $e_1\times\dots\times e_n:X_{g,n,d}\to X^n$ at the marked points.

\item Any polynomial of the first Chern classes $c^{(1)},\dots ,c^{(n)}$ of the
line bundles over $X_{g,n,d}$ consisting of tangent lines to the mapped curves
at the marked points. One defines these line bundles (by identifying the
Cartesian product of the forgetful and evaluation maps $X_{n+1,d}\to
X_{n,d}\times X$ with the {\em universal stable map} over $X_{n,d}$) as normal
line bundles to the $n$ embeddings $X_{n,d}\to X_{n+1,d}$ defined by the
$n$ marked points of the universal stable map. We will call these line
bundles {\em the universal tangent lines} at the marked points.

\item Pull-backs of cohomology classes of the Deligne - Mumford spaces
by contraction maps $\pi: X_{n,d}\to \calm _{0,n}$. We will make use of the
classes $A_I:= A_{i_1,...,i_k} $ Poincare-dual to fundamental cycles of fibers
of forgetful maps $\calm _{0,n}\to \calm _{0,k}$.

\end{enumerate}

We define the {\em GW-invariant}

\[ A_I\lan \phi _1,...,\phi _n \ran _{n,d} :=
\int _{X_{n,d}} \ \pi ^*A_I \w e_1^*\phi _1 \w ... \w e_n^*\phi _n .\]

It has the  following meaning in enumerative geometry: it counts the number
of pairs ``a degree-$d$ holomorphic map $\CC P^1 \to X$
with given $k$ points mapped to given $k$ cycles, a configuration of $n-k$
marked points mapped to the $n-k$ given cycles''.

\bigskip

Suppose now that the ample manifold $X$ is provided with a Killing action
of a compact Lie group $G$. Then $G$ act also on the moduli spaces of stable
maps. The evaluation, forgetful and contraction maps are $G$-equivariant,
and one can define correlators
$A_I\lan\phi_1,\dots \phi_n\ran_{n,d}$ of {\em equivariant}
cohomology classes of $X$.

The equivariant cohomology $H^*_G(M)$ of a $G$-space $M$ is defined as the
ordinary cohomology $H^*(M_G)$ of the homotopic quotient
$M_G=EG\times_G M$ --- the total space of the $M$-bundle
$p:M_G\to BG$ associated with the universal principal $G$-bundle
$EG\to BG$. The characteristic class algebra $H^*(BG)=H^*_G$(pt) plays the role
of the coefficient ring of the equivariant theory (so that $H^*_G(M)$ is a
$H^*_G$(pt)-module).  If $M$ is a compact manifold with smooth $G$-action, the
push-forward $p_*:H^*_G(M)\to H^*_G$(pt) (``fiberwise integration") provides
the equivariant cohomology of $M$ with intersection theory with values in
$H^*_G$(pt).

We introduce the {\em equivariant GW-invariants},
$A_I(\lan\phi_1,\dots ,\phi_n\ran_{n,d}$, with
values in $H^*(BG)$, where $\phi_1,\dots ,\phi_n\in H^*_G(X)$. Values of such
invariants on fundamental cycles of maps $B\to BG$ are accountable for
enumeration of rational holomorphic curves in families of complex manifolds
with
the fiber $X$ associated with the principal $G$-bundles over a
finite-dimensional
manifold $B$.

\section{The WDVV equation}
\label{sec3}
\setcounter{equation}{0}

One of the main structural results about Gromov-Witten invariants ---
the composition rule \cite{RT},\cite{MS} --- expresses all genus-0
correlators via the 3-fixed point ones, (we denote them $\lan \phi_1, ...,
\phi _n \ran _{n,d} $ since the corresponding $A_I=1$)
satisfying additionally the
so-called {\em Witten-Dijkgraaf-Verlinde-Verlinde} equation. We will see here
that the same result holds true for equivariant Gromov-Witten invariants
(at least in the ample case).

Following \cite{WD}, introduce the {\em potential}
\begin{equation}
\label{eq3.1}
F =\sum^\infty_{n=0} \ \frac{1}{n!} \ \sum_d q^d \
\lan t,\dots ,t\ran_{n,d} \ .
\end{equation}
It is a formal function of $t\in H^*_G(X)$ with values in the coefficient ring
$\L=H^*_G(pt,\CC[[q]]$). Here $\CC[[q]]$ stands for some  completion of
the group algebra $\CC[H_2(X,\ZZ)]$ so that the symbol
$q^d=q^{d_1}_1\dots q^{d_k}_k$ represents the class $(d_1,\dots ,d_k)$ in the
lattice $\ZZ^k=H_2(X,\ZZ)$ of 2-cycles. Fundamental classes of
holomorphic curves in $X$ have non-negative coordinates with respect to a basis
of Kahler forms so that the formal power series algebra $\CC[[q]]$ can be
taken on the role of the completion. Strictly speaking, the formula \ref{eq3.1}
defines $F$ up to a quadratic polynomial of $t$ since the spaces $X_{n,0}$ are
defined only for $n\geq 3$.

Denote $\nabla$ the gradient operator with respect to the equivariant
intersection pairing $\lan \ , \ \ran$ on $H^*_G(X)$. It is defined over the
field of fractions of $H^*_G(pt)$. The WDVV equation is an identity between
third directional derivatives of $F$. It says that

\begin{equation}
\label{eq3.2}
\lan\nabla F_{\a,\b},\nabla F_{\g,\d}\ran
\end{equation}

{\em is totally symmetric with respect to permutations of the four directions}
$\a,\b,\g,\d\in H^*_G(X)$.

\begin{thm} The WDVV equation holds for ample $X$.
\end{thm}

Notice that
\begin{equation}
\label{eq3.3}
\lan \nabla\int_X a\wedge t,\nabla\int_X b\wedge t\ran=
\lan a,b\ran
\end{equation}
has geometrical meaning of integration $\int_{\D\subset X\times X}
a\otimes b$ over the diagonal in $X\times X$.

In order to prove the non-equivariant version of the WDVV equation one
interprets the 4-point correlators $A_{1234}\lan\a,\b,\g,\d \ran_{4,d}$
which are totally symmetric in $\a,\b,\g,\d$
as integrals over the fibers $\G_\l$ of the contraction map
$\pi: X_{4,d}\to \calm _4 =\CC P^1$ and specializes the cross-ratio $\l$ to
$0, 1$ or $\infty$.
Stable maps corresponding to generic points of, say, $\G_0$ are glued from
a pair
of maps  $f_1: (\CC P^1,p_1,p_2,a_1)\to X$,
$f_2:(\CC P^1,p_3,p_4,a_2)\to X$ of degrees $d_1+d_2=d$ with
three marked points each, satisfying the diagonal condition $f_1(a_1)=
f_2(a_2)$. One can treat such a pair as a point in $X_{3,d_1}\times
X_{d_2,3}$ situated on the inverse image $\G_{d_1,d_2}$ of the diagonal
$\D \subset X\times X$ under the evaluation map $e_3\times e_3$.
The {\em glueing map} $\displaystyle{\sqcup_{d_1+d_2=d}} \G_{d_1,d_2}\to\G$ is
an isomorphism at generic points and therefore it identifies the analytic
fundamental cycles. This means that
\[
A_{1234}\lan\a,\b,\g\,d\ran_{4,d} =\sum_{d_1+d_2=d}
\lan\nabla\lan\a,\b,t \ran_{3,d_1},
    \nabla\lan\g,\d,t \ran_{3,d_2}\ran \ .
\]
The above argument applies to the correlators
$A_{1234}\lan\a,\b,\g\,d, t,\dots ,t\ran_{n+4,d}$ with additional marked points
and gives rise to
\begin{equation}
\label{eq3.4}
\lan\nabla F_{\a,\b},\nabla F_{\g,\d}\ran =\sum^\infty_{n=0} \
\frac{1}{n!} \ \sum_d \ q^d
A_{1234}\lan\a,\b,\g,\d,t\dots t\ran_{4+k,d}
\end{equation}
which is totally symmetric in $\a,\b,\g,\d$.

Ampleness of $X$ is used here only in order to make sure that the moduli
spaces
have fundamental cycles and that the diagonal in $X\times X$ consists of
regular values of the evaluation map $e_3\times e_3$.

In order to justify the above argument in the equivariant situation, it is
convenient to reduce the problem to the case of tori actions (using maximal
torus of $G$) and use the De Rham version of equivariant cohomology theory.

For a torus $G=(S^1)^r$ acting on a manifold $M$ the equivariant De Rham
complex
\cite{AB} consists of $G$-invariant differential forms on $M$ with
coefficients in
$\CC [u_1,\dots ,u_r]=H^*_G(pt)$, provided with the coboundary operator
$d_G=d+\sum^r_{s=1} u_si_s$ where $i_s$ are the operators of contraction
by the vector fields generating the action. Applying the ordinary Stokes
formula to $G$-invariant forms and $G$-{\em invariant} chains
we obtain   well-defined functionals $H^*_G(M)\to\CC [u]$ of
{\em integration over invariant cycles}. The identity \ref{eq3.4} follows
now from the obvious $G$-invariance of the analytic varieties $\Gamma_{\l }$,
$\Gamma$, $\Gamma_{d_1,d_2}$.

A similar argument proves a composition rule that reduces computation of all
equivariant correlators $A_I\lan ... \ran$ to that of $\lan ... \ran $.

\section{Ample vector bundles}
\label{sec4}

The following construction was designed by M. Kontsevich in order to extend the
domain of applications of WDVV theory to complete intersections in ample
Kahler manifolds.

Let $E\to X$ be an {\em ample} bundle, that is, a holomorphic vector
bundle spanned by its holomorphic sections. For stable
$f: (C,p)\to X$ (of degree $d$, with $n$ marked points), the spaces
$H^0(C,f^*E)$ form a holomorphic vector bundle $E_{n,d}$ over the moduli
space $X_{n,d}$.
If $f$ is glued from $f_1$ and $f_2$ as in the proof of
(\ref{eq3.4}), then $H^0(C,f^*E)= \ker (H^0(C_1,f^*_1E)
\oplus H^0(C_2,f^*_2E)
\stackrel{e_1-e_2}{\longrightarrow} e_1^*E=e_2^*E)$ where
$e_i: H^0(C_i,f^*_iE)\to e^*_iE$ is defined by evaluation of
sections at the marked point $a_i$.

This allows one to construct a solution $F$ to the WDVV equation
starting with an ample $G$-equivariant bundle $E$ and any invertible
$G$-equivariant multiplicative characteristic class $c$
(the total Chern class would be a good example).

Redefine
\[ \langle a,b\rangle \ := \ \int_X a\wedge b\wedge c(E) \ , \]
\[ \langle t,\dots ,t\rangle_{n,d} \ := \  \int_{X_{n,d}} e^*_1t \wedge
   \dots e^*_n t \wedge c(E_{n,d}) \ , \]
\[ F(t) \  = \  \sum^\infty_{n=0} \frac{1}{k!} \sum_d q^d
  \langle t,\dots ,t\rangle_{n,d} \ .\]
Then $\langle\nabla F_{\a,\b},\nabla F_{\g\d}\rangle$
{\em is totally symmetric in} $\a,\b,\g,\d$.

This construction bears a limit procedure from the total Chern class to the
(equivariant) Euler class, and the limit of $F$ corresponds to the GW-theory
on the submanifold $X'\subset X$ defined by an (equivariant) holomorphic
section $s$ of the bundle $E$. Namely, the section $s$ induces a
holomorphic section $s_{n,d}$ of $E_{n,d}$, and the (equivariant) Euler class
$Euler \ (E_{n,d})$  becomes represented by some cycle $[X'_{n,d}]$ situated
in the zero locus $X'_{n,d}:=s^{-1}_{n,d}(0)$ of the induced section.  The
variety $X'_{n,d}$ consists of stable maps to $X'$, the Euler cycle
$[X'_{n,d}]$
plays the role of the virtual fundamental cycle in $X'_{n,d}$, and the
correlators
\[ \lan t, ... ,t \ran _{n,d} := \int _{X_{n,d}} \ e_1^*t ... e_n^*t \ Euler \
(E_{n,d}) = \int _{[X'_{n,d}]} \ e_1^*t ... e_n^*t \]
are correlators of GW-theory on $X'$ between the classes $t$ which come from
the ambient space $X$.

\bigskip

Another solution of the $WDVV$-equation can be obtained from the bundles
$E'_{d,k}:= H^1(C,f^*E^*)$: one should put
$\langle a,b\rangle  := \int_X a\wedge b\wedge c^{-1}(E^*)$,
$\langle t,\dots ,t\rangle_{n,d}   = \int_{X_{n,d}}
  e^*_1t\wedge\dots\wedge e^*_n t\wedge c(E'_{n,d})$ for $d\neq 0$
and $\langle t,\dots ,t\rangle_{n,0}=\int_{X_{n,0}} e^*_1t\wedge
\dots\wedge e^*_n t\wedge c(E^*)$.

\section{Quantum cohomology}
\label{sec5}

One interprets the WDVV equation as the associativity identity for the
{\em quantum cup-product} on $H^*_G(X)$ defined by
\[
\langle \a*\b,\g\rangle = F_{\a,\b,\g} \ .
\]
It is a deformation of the ordinary cup-product
(with $t$ and $q$ in the role of parameters)  in the category
of (skew)-commutative algebras {\em with unity}:
\begin{equation}
\label{eq5.1}
\langle \a*1,\g\rangle = \langle\a,\g\rangle \ .
\end{equation}
Indeed, the push-forward by the forgetful map
$\pi : X_{n,d}\to X_{n-1,d}$ (with $n\geq 3$) sends
$1\in H^*_G(X_{n,d})$ to 0 unless $d=0$ and $k=3$ in which case
$X_{n,d}=X$ and $X_{n-1,d}$ is not defined.

The structure usually referred in the literature as the
{\em quantum cohomology algebra} corresponds to the restriction of the
deformation $*_{t,q}$ to $t=0$. As it is shown in \cite{KM}, in many
cases the function $F$ can be recovered on the basis of WDVV-equation from
the structural constants
$F_{\a,\b,\g}|_{t=0}(q)$ of the quantum cohomology algebra due to the following
symmetry of the potential $F$. Let $u \in H^2_G(X)$ and $(u_1,...,u_k)$ be
its coordinates with respect to the basis of the lattice
$(\ZZ ^k)^*=H^2(X)=H^2_G(X)/H^2_G(pt)$ (so that $u_i \in H^*_G(pt)$).
Then
\begin{equation}
\label{eq5.2}
(F_{\a , \b ,\g })_u=\sum_{i=1}^k u_iq_i \p F_{\a , \b ,\g }/\p q_i
\quad \forall \a,\b,\g \in H^*_G(X) \ .
\end{equation}
The identity (\ref{eq5.2}) follows from the obvious push-forward formula
$\pi_*u=d_i u_i$.

The symmetry $(6)$ can be interpreted in the way that the quantum deformation
of the cup-product restricted to $t=0$ is equivalent to the deformation with
$q=1$ and $t$ restricted to the $2$-nd cohomology of $X$ (in the equivariant
setting it is better however to keep both parameters in place ---
see Sections $7, 8$).

\bigskip

In this paper, we will use the term {\em quantum cup-product} for the entire
$(q,t)$-deformation and {\em reserve the name {\em quantum cohomology algebra}
for the restriction of the quantum cup-product to $t=0$}.

I have heard some complaints about such terminology because it allows many
authors to compute quantum cohomology algebras without even mentioning the
deformation in $t$-directions. There are some indications however that (despite
the equivalence $(6)$) the $q$-deformation has a somewhat different nature than
the $t$-deformation. The loop space approach \cite{HG1} and our computations
in Sections $9$ -- $11$ seem to emphasize this distinction.

\bigskip

Quantum cohomology algebras of the classical flag manifolds have been
computed in
\cite{GK}, \cite{K} on the basis of several
conjectures about properties of $U_n$-{\em equivariant}
quantum cohomology (see also \cite{AS} where a slightly different formalism
was applied). The answer (in terms of generators and relations)
for complete flag manifolds  $U_n/T^n$ is strikingly related to
conservation laws of Toda lattices. The conjectures named
in \cite{GK} the {\em product, induction} and {\em restriction}
properties and describing behavior of equivariant quantum cohomology  under
some natural constructions, were motivated by interpretation of the
quantum cohomology in terms of Floer theory on the loop space $LX$.
Although a construction of the equivariant counterpart of the Floer - Morse
theory on $LX$ remains an open problem, the three conjectured axioms can
be justified within the Gromov-Witten theory. This was done by B.Kim \cite{K2}.
The induction and
restriction properties follow directly from definitions given in this paper
and hold for the entire quantum deformation (not only at $t = 0$), while the
``product'' axiom that the $G_1 \times G_2$-equivariant quantum cohomology
algebra of $X_1\times X_2$ is the tensor product of the $G_i$-equivariant
quantum cohomology algebras of the factors $X_i$  has been verified in
\cite{K2} for ample manifolds.

Behavior of the quantum cup-product at
$t\neq 0$ under the Cartesian product operation on the target manifolds is
much more complicated than the operation of the tensor product.

\section{Floer theory and $D$-modules}
\label{sec6}

Structural constants $\langle \a * \b,\g\rangle$ of the quantum cup-product
are derivatives $\partial_{\b}F_{\a,\g}$ of the same function.  This allows to
interpret the WDVV-equation as integrability condition of some connections
$\nabla_\hbar $ on the tangent bundle $T_H$ of the space $H = H^*(X,{\Bbb C})$.
Namely, put $t = \sum t_{\a}p_{\a}$ where $p_1 = 1,p_2,\dots,p_N$ is a basis
in $H$ and define
\[
\nabla_\hbar = \hbar d - \sum(p_{\a}*)dt_{\a} \wedge:
\Omega^0(T_H) \to \Omega^1(T_H)
\]
where $p_{\a}*$ are operators of quantum multiplication by $p_{\a}$.  Then
$\nabla_\hbar \circ \nabla_\hbar = 0$ {\em for each value of the parameter}
$\hbar$.
Notice that the integrability condition that reads ``the system of differential
equations $\hbar \partial_{\a}I = p_{\a} * I$ has solutions $I \in
\Omega^0(TH)$'' is actually obtained as a somewhat combinatorial statement
(the WDVV-equation) about coefficients of the series $F$.

In \cite{HG1}, \cite{HG} we attempted to improve this unsatisfactory
explanation of the integrability property by
describing a direct geometrical meaning of the solutions $I$ in terms of
$S^1$-equivariant Floer theory on the loop space $LX$.  Briefly, the universal
covering $\widetilde{LX}$ carries the action of the covering transformation
lattice $\pi_2(X)$ with generators $q_1,\dots,q_k$ and the $S^1$-action by
rotation of loops which preserves natural symplectic forms $\o_1,\dots,\o_k$ on
$LX$ and thus defines corresponding Hamiltonians $H_1,\dots,H_k$ on
$\widetilde{LX}$ (the action functionals).  The Duistermaat--Heckman forms
$w_i+ \hbar H_i$ (here $\hbar $ is the generator of $H_{S^1}^*(pt)$)
are equivariantly closed,
and operators $p_i$ of exterior multiplication by these forms have the
following Heisenberg commutation relations with the covering transformations:
\[
p_iq_j - q_jp_i = \hbar q_j\d_{ij}.
\]
Conjecturally, this provides $S^1$-equivariant Floer cohomology of
$\widetilde{LX}$ with a ${\cald}$-module structure which is equivalent to the
above system of differential equations (restricted to $t = 0$, $q \ne 0$) and
reduces to the quantum cohomology algebra in the quasi-classical limit
$\hbar = 0$ (see \cite{HG1, GK}).

In this section we describe solutions to $\nabla_\hbar I = 0$ by imitating the
$S^1$-equivariant Floer theory (which is still to be constructed) within the
framework of Gromov--Witten theory. This construction turns out to be crucial
in our proof in Section 11 of the mirror conjecture for Calabi-Yau projective
complete intersections.

\bigskip

One may think of the graph of an algebraic
loop ${\Bbb C}P^1{\backslash}\{0,\infty\} \to X$ of degree $d$ as of a stable
map ${\Bbb C}P^1 \to X \times {\Bbb C}P^1$ of bidegree $(d,1)$.  Our starting
point consists in interpretation of the moduli space $L_d(X)$ of such stable
maps as a degree-$d$ approximation to $\widetilde{LX}$ and application of
equivariant Gromov--Witten theory to the action of $S^1$ on the second factor
${\Bbb C}P^1$ with the fixed points $\{0,\infty\}$.

In the theorem below we assume $X$ to be ample.  It is natural to expect
however that the theorem holds true whenever the non-equivariant Gromov--Witten
theory works for $X$ since the $S^1$-action is non-trivial only on the factor
${\Bbb C}P^1$ which is ample on its own.

Let $\langle \ ,\ \rangle$ be the Poincare pairing on $H = H^*(X,{\Bbb C})$.
The equivariant cohomology algebra $H_{S^1}^*(X \times {\Bbb C}P^1)$ is
isomorphic to $H \otimes_{\Bbb C} {\Bbb C}[p,\hbar ]/(p(p-\hbar ))$ with the
$S^1$-equivariant pairing
\[
(\varphi,\psi) = \frac {1}{2\pi i} \oint \frac {\langle \varphi,\psi\rangle
dp}{p(p-\hbar )}.
\]
Localization in $\hbar $ allows to introduce coordinates $\varphi = tp/\hbar +
\tau(\hbar -p)/\hbar $, $\tau,t \in H$, diagonalizing the equivariant pairing:
\[
((\tau,t),(\tau ',t')) = \frac {\langle t,t' \rangle  - \langle
\tau,\tau '\rangle }{\hbar }.
\]
The potential 
${\calf}(t,\tau,h,q,q_0)$ satisfying the
equivariant WDVV-equation for $X \times {\Bbb C}P^1$ expands as
\[
{\calf} = {\calf}^{(0)} + q_0{\calf}^{(1)} + q_0^2{\calf}^{(2)}\dots
\]
according to contributions of stable maps of degree $0,1,2,\dots$ with respect
to the second factor.  Denote $F = F(t,q)$ the potential (\ref{eq3.1}) of
the $GW$-theory for $X$.

\begin{thm}
{\em (a)} ${\calf}^{(0)} = (F(t,q)-F(\tau,q))/\hbar $.

{\em (b)} The matrix $(\Phi_{\a\b}) :=
(\partial^2{\calf}^{(1)}/\partial\tau_{\a}\partial t_{\b})$ is a fundamental
solution of $\nabla_{\pm \hbar }I = 0$:
\[
-\hbar \frac {\partial}{\partial \tau_{\g}} \Phi = p_{\g}(t)\Phi \ ,\]
\[ \hbar \frac {\partial}{\partial t_{\g}} \Phi^* = p_{\g}(\tau)\Phi^* \ ,\]
where $p_{\g} = (p_{\a}^{\b})_{\g}$, $\g = 1,\dots,N$, are matrices of quantum
multiplication by $p_1 = 1,\dots,p_N$, and $\Phi^*$ is transposed to $\Phi$.
\end{thm}

\pf Moduli spaces of bidegree-$(d,0)$ stable maps to $X \times {\Bbb C}P^1$
coincide with $X_{n,d} \times {\Bbb C}P^1$.  This implies (a) and shows that
the WDVV-equation for $\calf$ {\em modulo} $q_0$ follows from the WDVV-equation
for $F$.  Part (b) follows now directly from the WDVV-equation for $\calf$
{\em modulo} $q_0^2$ and from
\[
\hbar \frac{\p }{\p t_1} \Phi _{ \a \b } = \Phi _{\a \b } =
-\hbar \frac{\partial }{\partial \tau _1 } \Phi _{\a \b }
\]
due to (\ref{eq5.1}) and (\ref{eq5.2}). Here $\p /\p t_1, \p/\p \t _1$ are
derivatives in the direction $1\in H^*(X)$ of the identity components of
$t$ and $\t$ respectively.

\bigskip

The following corollary is obtained by expressing equivariant correlators
$\Phi_{\a\b}$ via localization of equivariant cohomology classes of moduli
spaces $L_d(X)$ to fixed points of the $S^1$ action.

Define
\begin{equation}
\label{eq6.1}
\psi_{\a\b} = \sum_{n=0}^{\infty} \frac {1}{n!} \sum_d q^d\langle \frac
{p_{\b}}{\hbar +c} ,t,\dots,t,p_{\a}\rangle_{n+2,d}
\end{equation}
where $c$ is the first Chern class of the line bundle over $X_{k,d}$ introduced
in Section $1$ as ``the universal tangent line at the first marked
point'', and $\langle \frac {p_{\b}}{\hbar +c}, p_{\a}\rangle_{2,0} :=
\langle p_{\a}, p_{\b}\rangle $.

\begin{cor}
$\hbar \p \psi /\p t_{\g } = p_{\g}(t)\psi$,
i.e., the matrix $\psi$ is (another) fundamental solution of
$\nabla_{\hbar } I = 0$.
\end{cor}

\pf A fixed point in $L_d(X)$ is represented by a stable map $C_0 \cup {\Bbb
C}P^1 \cup C_{\infty} \to X \times {\Bbb C}P^1$ where $\varphi_i: C_l \to X
\times \{i\}$ are stable maps of degrees $d_1 + d_2 = d$ connected by a
``constant loop'' ${\Bbb C}P^1 \stackrel{\simeq}{\rightarrow} \{x\} \times
{\Bbb C}P^1$.  Thus components of $L_d(X)^{S^1}$ can be identified with
submanifolds in $X_{d_1,k_1+1}^{(0)} \times X_{d_2,k_2+1}^{(\infty)}$ defined
by the diagonal constraint $e_1(\varphi_0) = e_1(\varphi_{\infty})$, with
$\hbar^2(\hbar + c(0))(\hbar - c(\infty))$ to be the
equivariant Euler class of the normal bundle.  This gives rise to
\begin{equation}
\label{eq6.2}
\hbar^2\Phi_{\a\b} = \sum_{\e,\e '}
\psi_{\a\e}(\tau ,-\hbar )\ \eta^{\e\e'} \psi_{\e '\b}(t,h)
\end{equation}
where $\sum \eta^{\e\e'}p_{\e} \otimes p_{\e'}$ is the coordinate expression
of the diagonal cohomology class of $X\subset X \times X$.

\bigskip

We give here several reformulations which will be convenient for computation
of quantum cohomology algebras in Sections $9$ -- $10$.

Consider the specialization of the connection $\nabla _{\hbar }$ to the
parameter subspace corresponding to the deformation of the
quantum cup-product along the $2$-nd cohomology  (this is accomplished by
putting first $t=0$ and then replacing $q^d$ by $\exp (d,t)$ where
$t=(t_1,...,t_k)$
represents coordinates on $H^2(X)$ with respect to the basis
$p^{(1)},..., p^{k)} \in H^2(X)$ . In this new setting put
\[ s_{\a, \b} := \sum _d e^{dt} \lan p_{\b } \frac{e^{pt/\hbar }}{\hbar + c},
p_{\a }\ran \]
where $pt:= \sum p^{(i)}t_i$.

{\bf Corollary 6.3.} {\em The matrix $(s_{\a ,\b }(t))$ is a fundamental
solution to
\[ \nabla _{\hbar } \ s = 0 : \ \hbar \frac{\p}{\p t_i} \ s = p^{(i)} * s .\] }

{\em Proof.} One should combine Corollary $6.2$ with iterative applications
of the following generalized symmetries $(5),(6)$:
\[ \lan f(c), ..., 1 \ran _{n+1, d} = \lan \frac{f(0)-f(c)}{c},... \ran _{n,d}
, \]
\[ \lan f(c), ..., p^{(i)} \ran _{n+1, d} = d_i \lan f(c),... \ran _{n,d} +
\lan p^{(i)} \frac{f(0)-f(c)}{c} , ... \ran _{n,d} .\]
Here $f(c)$ is a function of $c$ with values in $H^*(X)$.

The symmetries are easily verified on the basis of the following geometrical
properties of universal tangent lines:

(i) Consider the push-forward along the map $\pi : X_{n+1,d}\to
X_{n,d} $ (forgetting the last marked point). It is easy to see that
the difference $\pi ^*(c)-c$ between the Chern class of the universal
tangent line at the $1$-st marked point and the pull-back of its counterpart
from $X_{n,d}$ is represented by the fundamental cycle of the section
$i: X_{n,d}\to X_{n+1,d}$ defined by the first marked point.

(ii) $i^* (c) = c$.

In particular $\pi _* (1/(\hbar +c)) = 1/[\hbar (\hbar +c)] $.

{\bf Corollary 6.4.} {\em Consider the functions
\[ s_{\b }:= \sum _d e^{dt }
\lan p_{\b } \frac{e^{pt/\hbar }}{\hbar + c}, 1\ran _{2,d} \ .\]
Let $P(\hbar \p/\p t, \exp t, \hbar)$ be a differential operator
annihilating simultaneously all the functions $s_{\a }$. Then
the relation $P(p^{(1)},..., p^{(k)} ,q_1,...,q_k , 0)=0$ holds
in the quantum cohomology algebra of $X$ (we assume here that $P$
depends only on non-negative powers of $\hbar $).}

The functions $s_{\b }$ --- ``the first components of the vector-solutions
$s_{\a, \b}$'' --- generate a left $\cal{D}$-module with the {\em solution}
locally constant sheaf described by the flat connection $\nabla _{\hbar }$ and
with the characteristic variety isomorphic to the spectrum of the quantum
cohomology algebra.

\bigskip

All results of this section extend literally to the equivariant setting and/or
to the generalization to ample vector bundles described in Section $4$. We will
apply them in this extended form in Sections $9$--$11$.

\bigskip

{\em Remarks.}
$ 1$) The universal formula (\ref{eq6.1}) for solutions of $\nabla_\hbar I = 0$
was perhaps discovered independently by several authors. I first learned this
formula from  R. Dijkgraaf. It can also be found in \cite{Db} in the
{\em axiomatic} context of conformal topological field theory. One can prove
it directly from a recursion relation (in the spirit of WDVV-equation) for
so-called {\em gravitational descendents} --- correlators involving the
first Chern classes of the universal tangent lines (or, in a slightly
different manner, by describing explicitly the divisor in $X_{n,d}$
representing $c$).
Our approach provides an interpretation of (\ref{eq6.1}) in
terms of fixed point localization in equivariant cohomology.

$ 2$)  One can generalize our theorem to bundles over ${\Bbb C}P^1$ with the
fiber $X$.  This seems to indicate that a straightforward ``open-string''
approach to $S^1$-equivariant Floer theory on $\widetilde{LX}$ would be more
powerful and flexible than the approximation by Gromov--Witten theory on $X
\times {\Bbb C}P^1$ described above.

$3$) Although the theorem provides a geometrical interpretation of
solutions to $\nabla_\hbar I = 0$,
it does not eliminate the combinatorial
nature of the integrability condition.  Indeed, the theorem is deduced from an
equivariant WDVV-equation which in its turn can be interpreted as an
integrability condition.  Of course one can explain it using the $S^1 \times
S^1$-equivariant WDVV-equation on $(X \times {\Bbb C}P^1) \times {\Bbb C}P^1$,
etc.  It would be interesting to find out whether this process converges.

\section{Frobenius structures}
\label{sec7}

In \cite{Db}, B. Dubrovin studied geometrical structures defined by
solutions of WDVV-equations on the parameter space and reduced classification
of
generic solutions to the classification of trajectories of some Euler-like
non-autonomous Hamiltonian systems on $so_N^*$.  We show here how this approach
to equivariant Gromov--Witten theory yields analogous Hamiltonian systems on
the
affine Lie coalgebras $\widehat{so}_N^*$.

The quantum cup-product on $H = H_G^*(X)$ considered as an $N$-dimensional
vector space over the field of fractions $K$ of the algebra $H_G^*(pt)$ defines
a formal {\em Frobenius structure} on $H$.  The structure consists of the
following ingredients.

\begin{enumerate}
\item A symmetric $K$-bilinear inner product $\langle\ ,\ \rangle$,

\item a (formal) function $F: H \to K$ whose third directional derivatives
$\langle
a*b,c\rangle := F_{a,b,c}$ provide tangent spaces $T_tH$ with the Frobenius
algebra structure (i.e. associative commutative multiplication $*$
satisfying $\lan a*b, c\ran = \lan a, b*c \ran $).

\item The constant vector field $1\!\!1$ of unities of the algebras
$(T_tH, *)$ whose flow
preserves the multiplication $*$ (i.e. $L_{1\!\!1} (*)=0$).

\item Grading: In the non-equivariant case axiomatically studied
by B. Dubrovin it can be
described by the {\em Euler} vector field $E$, such that
the tensor fields $1\!\!1$, $*$ and $\lan \ , \ \ran$ are homogeneous (i.e.
are eigen-vectors of the Lie derivative $L_E$) of degrees $-1, 1$ and $D$
respectively (where $D = \dim _{\Bbb C} X$ in the models arising from
the GW-theory). In the equivariant GW-theory this grading axiom should be
slightly modified since the grading of the
structural ring $H_G^*(pt)$ is non-trivial and thus the natural Euler
operator $L_E$ is $\CC$-linear but not $K$-linear.
\end{enumerate}

\bigskip

The fact that the multiplication $*$ is defined on tangent vectors to $H$
means that the algebra $(\Omega ^0(T_H), *)$ can be naturally considered as
the algebra $K[L]$ of regular functions on some subvariety $L\subset T^*H$
in the cotangent bundle.
A point $t\in H$ is called {\em
semi-simple} if the algebra $(T_{t}H,*)$ is semi-simple, that is if $L \cap
T_{t}^*H$ consists of $N$ linearly independent points.

Flatness of the connection (defined on $T_H$)
\begin{equation}
\label{eq7.1}
\nabla_\hbar = \hbar d - \sum_{\a} p_{\a} * dt_{\a}
\end{equation}
implies \cite{GK} that $L$ is a Lagrangian submanifold in $T^*H$ near a
semisimple $t$.  Following \cite{Db}, introduce local {\em canonical
coordinates} $(u_1,\dots,u_N)$ such that the sections $(du_1,\dots,du_N)$ of
$T^*H$ are the $N$ branches of $L$ near $t$, and transform the connections
$\nabla_\hbar$ to these local coordinates and to a (suitably normalized)
basis $f_1,\dots,f_N$ of vector field on $H$ diagonalizing the $*$-product.

The result of this transformation can be described as follows.

(a) The basis $\{f_i\}$ can be normalized in a way that in the transformed
form
\begin{equation}
\label{eq7.2}
\nabla _{\hbar } = \hbar d - \hbar A^1 \wedge - D^1\wedge
\end{equation}
of the connection $\nabla _{\hbar} $ with
$D^1 = \operatorname{diag}(du_1,\dots,du_N)$, and $A_{ij} =
V_{ij}(u) d(u_i-u_j)/(u_i-u_j)$ for all $i \ne j$,
we will have additionally $A_{ii} =0 \ \forall i$.

(b) The vector field $1\!\!1$ in the canonical coordinates assumes the form
$\sum_k \partial_k$ where $\partial_k := \partial/\partial u_k$ are
the canonical idempotents of the $*$-product:
\begin{equation}
\label{eq7.3}
\partial_i * \partial_j = \d_{ij}\partial_j.
\end{equation}

(c) The (remaining part of the) integrability condition $\nabla_{\hbar} ^2 =0$
reads
$d(A^1) = A^1 \wedge A^1$ or
\begin{equation}
\label{eq7.4}
\partial_i\phi_{\a}^j = \phi_{\a}^i V_{ij}/(u_i-u_j),\ i \ne j,
\end{equation}
where $(\phi_{\a}^j)$ is the transition matrix, $\partial/\partial t_{\a} =
\sum_i \phi_{\a}^i f_i$; it can be reformulated as compatibility of the
PDE system (\ref{eq7.4}) for $(\phi_{\a}^j)$ completed by
\begin{equation}
\label{eq7.5}
\sum_k \partial_k\phi_{\a}^j = 0.
\end{equation}

(d) The Frobenius property $\langle a*b,c\rangle = \langle a,b*c\rangle$ of the
$*$-product shows that the diagonalizing basis $\{f_i\}$ is orthogonal, that
its normalization by $\langle f_i,f_j \rangle = \d_{ij}$ obeys $A_{ii} = 0$
and, additionally, implies anti-symmetricity $A_{ij} = -A_{ji}$, or
\begin{equation}
\label{eq7.6}
V_{ij} = -V_{ji}.
\end{equation}

The presence of the grading axiom (4) of Frobenius structures over $K = {\Bbb
C}$ allows B.Dubrovin to describe anti-symmetric matrices
$V = (V_{ij}) \in so_N^*$
satisfying the integrability conditions (\ref{eq7.4}) and (\ref{eq7.5}) in
{\em quasi-homogeneous} canonical coordinates (i.e. $L_E u_i = u_i$ so that
$E = \sum u_k\partial_k$)
as trajectories of $N$ commuting non-autonomous Hamiltonian
systems (see \cite{Db}):
\[
\partial_iV = \{H_i,V\}
\]
where the Poisson-commuting non-autonomous quadratic Hamiltonians $H_i$
on $so_N^*$ are given by
\[
H_i = \frac {1}{2} \sum_{j \ne i} \frac {V_{ij}V_{ji}}{u_i-u_j} .
\]

Consider now the following model modification of the grading axiom:
$K = {\Bbb C}[[\lambda^{\pm 1}]]$, $\operatorname{deg} \lambda = 1$.  In
quasi-homogeneous canonical coordinates $(u_1,\dots,u_N,\lambda)$ the Euler
vector field takes then on
\begin{equation}
\label{eq7.7}
L_E = \sum_k u_k\partial_k + \lambda\partial_{\lambda}.
\end{equation}

Introduce the connection operator
\[
{\Bbb V} = \lambda\partial_{\lambda} - V \in \widehat{so}_N^*
\]
and the qudratic Hamiltonians on the Poisson manifold $\widehat{so}_N^*$
\begin{equation}
\label{eq7.8}
{\calh}_i({\Bbb V}) = \oint H_i(V) \frac {d\lambda}{\lambda} .
\end{equation}

\begin{prop}
The Hamiltonians ${\calh}_1,\dots,{\calh}_N$ are in involution.  The operator
${\Bbb V}$ of a Frobenius manifold over $K$ satisfies the non-autonomous system
of Hamiltonian equations
\begin{equation}
\label{eq7.9}
\partial_i{\Bbb V} = \{{\calh}_i,{\Bbb V}\},\ i = 1,\dots,N.
\end{equation}
The columns $\phi_{\a} = (\phi_{\a}^i)$ of the transition matrix are
eigen-functions of the connection operator ${\Bbb V}$:
\begin{equation}
\label{eq7.10}
{\Bbb V}\phi_{\a} := (\lambda\partial_{\lambda} - V)\phi_{\a} = \left( \frac
{n}{2} - \operatorname{deg} t_{\a} + 1\right) \phi_{\a}.
\end{equation}
\end{prop}

\pf It can be obtained by a straightforward calculation quite analogous
to that in \cite{Db}.

\bigskip

In our real life the model equations (\ref{eq7.7}--\ref{eq7.10}) describe the
structure of Frobenius manifolds {\em over each semi-simple orbit of the
grading Euler field in the ground parameter space.  This parameter space is the
spectrum of the coefficient algebra} $H_G^*(pt,{\Bbb C}) \otimes {\Bbb
C}[q^{\pm 1}_1,\dots,q^{\pm 1}_k]$ (its field of fractions can be taken on the
role of the ground field $K$).  An orbit of the Euler vector field in this
parameter space is
semi-simple if the corresponding ${\Bbb C}$-Frobenius algebras are
semi-simple.

The equations (\ref{eq7.7}--\ref{eq7.10}) over semi-simple Euler orbits should
be complemented by the additional symmetries (\ref{eq5.2}).

In the next section we will show how the canonical coordinates of the axiomatic
theory of Frobenius structures emerge from localization formulas in equivariant
Gromov--Witten theory.

\section{Fixed point localization}
\label{sec8}

We consider here the case of a circle $T^1$ acting by Killing transformations
on a compact Kahler manifold $X$ with {\em isolated} fixed points only.
The case of tori actions with isolated fixed points requires only slight
modification of notations which we leave to the reader. Our
results are rigorous for ample $X$ (which includes homogeneous Kahher spaces of
compact Lie group and their maximal tori) while applications to general toric
manifolds (which are typically not ample) yet to be justified.

It is the Borel localization theorem that reduces computations in
torus-equivariant
cohomology to computations near fixed points. Let $\{ p_{\a}\} ,
\a=1,...,N$, be the fixed points of the action. We will denote with the same
symbols $p_{\a}$ the equivariant cohomology class of $X$ which restricts
to $1\in H^*_T(p_{\a})$ at $p_{\a}$ and to $0$ at all the other fixed points.
These classes
are well-defined over the field of fractions $\CC (\l )$
of the coefficient ring $H^*_T(pt)=\CC [\l ]$ and form the basis
of canonical idempotents in the semi-simple algebra $H^*_T(X,\CC (\l )$.
The equivariant Poincare pairing reduces to $\lan p_{\a},p_{\b}\ran =
\d _{\a,\b}/e_{\b}$ where $e_{\a}\in \CC [\l]$ is the equivariant Euler class
of the normal ``bundle'' $T_{p_{\a}}X\to p_{\a}$ to the fixed point.

The results described below apply to the setting of Section $4$ of a manifold
$X$ provided with an ample vector bundle in which case $e_{\a}$'s should be
modified accordingly.

\bigskip

The same localization theorem reduces computation of GW-invariants to that
near the fixed point set (orbifold) in the moduli spaces $X_{n,d}$.
A fixed point in
the moduli space is represented by a stable map to $X$ of a (typically
reducible) curve $C$ such that each component of $C$ is mapped to (the
closure of) an orbit of the complexified action $T_{\CC }:X$. Any such an
orbit is either one of the fixed points $p_{\a}$ or isomorphic to $(\CC-0)$
connecting two distinct fixed points $p_{\a}$ and $p_{\b}$ corresponding to
$0$ and $\infty $. Respectively, there are two types of components of $C$:

\noindent (i) Each component of $C$ which carries $3$ or more special points
must be mapped to one of the fixed points $p_{\a}$.

\noindent (ii) All other components are multiple covers $z\mapsto z^d$ of the
non-constant orbits, and their special points may correspond only to $z=0$ or
$\infty$.

The {\em combinatorial structure}
of such a stable map can be described by a tree
whose edges correspond to {\em chains} of components of type (ii) and should
be labeled by the total degree of this chain as a curve in $X$, and vertices
correspond to the ends of the chains. The ends may carry $0$ or $1$ marked
point, or correspond to a (tree of) type-(i) components with $1$ or more
marked points and should be labeled by the indices of these marked points
and by the target point $p_{\a}$.

{\em The fixed stable maps with different combinatorial structure belong to
different connected components of the fixed point orbifold in $X_{n,d}$.}

The results below are based on the observation that a stable map with the
first $k\geq 3$ marked points in a {\em given} generic configuration
(i.e. with the given generic value of the contraction map $X_{n.d}\to
\calm _{0,k}$) must have in an irreducible component $C_0$ in the underlying
curve $C$ which contains this given configuration of $k$ {\em special} points,
(so that the corresponding first $k$ marked points are located on the branches
outgoing these special points of $C_0$). The cause is hidden in the definition
of the contraction map (see \cite{Kn, BM}).

We will call the component $C_0$ {\em special}.

The observation applied to a fixed stable map of the circle action allows to
subdivide all fixed point components in $X_{3+n,d}$ into $N$ {\em types}
$p_i$ according
to the fixed points $p_i$ where the special component is mapped to. We
introduce the superscript notation $(...)^i$ for the contribution (via
Borel's localization formulas) of type-$p_i$ components into various
equivariant correlators. For example,
\[ F^i_{\a \b \g } =\sum _{n}\frac{1}{n!} \sum_d q^d
\lan p_{\a },p_{\b },p_{\g },t,...,t\ran _{3+n,d}^i \]
where $t=\sum _{\a=1}^N t_{\a}p_{\a}$ is the general class in
$H^*_T(X,\CC(\l))$, so that $F_{\a\b\g}=\sum_i F_{\a\b\g}^i$.

We introduce also the notations
\begin{itemize}
\item $\Psi _{\a\b}^i$ --- for contributions to $e_iF_{\a\b i}^i$ of those
fixed points which have the third marked point situated directly on the
special component $C_0$ (it is convenient here to introduce the normalizing
factor $e_i\in H^*_T(pt)$, the Euler class of the normal ``bundle'' to the
fixed point $p_i$ in $X$);
\item $\Psi _{\a}^i := \Psi _{\a 1\!\!1}^i=\sum_\b \Psi _{\a\b}^i$;
\item $D_{\a}^i$ --- for contributions to $e_iF_{\a i i}$ of those
fixed points which have the second and third marked points situated directly
on the special component;
\item $\D ^i$ --- for contributions to $e_iF_{i i i}$ of those fixed
points which have the first three marked points situated directly on the
special component;
\item $u_i=t_i+$ contributions to $e_iF_{i i}$ of all those fixed point
components in $X_{2+n,d}$ for which the first two marked points belong to
the same vertex of the tree describing the combinatorial structure.
\end{itemize}

The correlators $u_i$ can be also interpreted as contributions to the
{\em genus-$1$} equivariant correlators
\[ \sum _n\frac{1}{n!}\sum_d q^d (t,...,t)_{n,d} \]
with {\em given} complex structure of the elliptic curve of those $T$-invariant
classes which map the (only) genus-$1$ component of the curve $C$ to the
fixed point $p_i$.

\bigskip
{\bf Theorem 8.1.} {\em (a) The functions $u_1(t),...,u_N(t)$ are the canonical
coordinates of the Frobenius structure on $H^*_G(X,\CC (\l))$.

(b) The functions $D_{\a}^i(t)$ are eigen-values of the quantum multiplication
by $p_{\a}$: $du_i=\sum _{\a} D_{\a }^i dt_{\a} $.

(c) The transition matrix $(\Psi_{\a}^i)$ provides simultaneous diagonalization
of the quantum cup product: $F_{\a\b\g}^i=\Psi_{\a}^i D_{\b}^i \Psi_{\g}^i$
and obeys the following orthogonality relations:
\[ \sum_i \Psi_{\a}^i\Psi_{\b}^i=\d _{\a\b}/e_{\b},
\sum _{\a} \Psi_{\a}^i\Psi_{\a}^j =\d _{ij} \ .\]

(d) The Euclidean structure on the cotangent bundle of the Frobenius
manifold (defined by the equivariant intersection pairing in $H^*_T(X)$)
in the canonical coordinates $u_i$ takes on $\lan du_i,du_j \ran =
(\D^i)^2\d_{ij}e_j$ and additionally
\[ (\D ^i)^{-1}=\sum_{\a}\Psi_{\a}^i, \ \Psi_{\a}^i=\frac{D_{\a}^i}{\D^i},
\ \Psi_{\a\b}^i=\frac{D_{\a}^i D_{\b}^i}{\D^i} \ .\]
  }

\bigskip

{\em Proof.} We first apply the localization
formula
\[ A_{1234}\lan ... \ran _{n,d}=\sum _i  A_{1234}\lan ...\ran _{n,d}^i \]
to the $4$-point equivariant correlators with the fixed cross-ratio $z$
of the $4$ marked points and only after this specialize the cross-ratio
to $0$,$1$ or $\infty$. This gives rise to the {\em local} WDVV-identities
\[ \Psi_{\a\b}^i\Psi_{\g\d}^i \ \text{\em is totally symmetric in} \ \
{\a,\b,\g,\d} \]
which is independent of the global WDVV-equation.
When combined with the global identities
\[ A_{1234}\lan 1\!\!1,p_{\a},p_{\b},p_{\g} \ran _{n,d}
=\lan p_{\a},p_{\b},p_{\g} \ran _{n,d} \]
they yield the orthogonality relation $\sum_i \Psi _{\a}^i\Psi_{\b}^i=
\d_{\a\b}/e_{\b}$ and localization formulas
\[ F_{\a\b\g}=\sum_i\Psi_{\a\b}^i\Psi_{\g}^i \]
for the structural constants of the quantum cup-product.

A similar argument with $>4$-point correlators $A_{12345...}\lan ...\ran ^i$
proves the diagonalization
\[ \lan p_{\a}*p_{\b},p_{\g}\ran =
\sum_i \Psi_{\a}^iD_{\b}^i\Psi_{\g}^i/e_i  \ ,\]
\[ \lan p_{\a}*p_{\b}*p_{\g},p_{\d}\ran =
\sum _i \Psi_{\a}^i D_{\b}^i D_{\g}^i \Psi_{\d}^i/e_i \ \]
and the identities
\[ \Psi_{\a\b}^i=\Psi_{\a}^iD_{\b}^i, \
(\D ^i)^{-1}=\sum_{\a } \Psi_{\a}^i \ .\]

Finally, the identity $du_i=\sum_{\a} D_{\a}^idt_{\a} $ follows
directly from the definition of $u_i$ and implies that $u_1,...,u_N$
are the canonical coordinates of the Frobenius structure.

\section{ Projective complete intersections}
\label{sec9}

We are going to describe explicitly solutions of the
differential equations arising from quantum cohomology of projective complete
intersections. Lex $X$ be such a non-singular complete intersection in
$Y:=\CC P^n $ given by $r$ equations of the degrees $(l_1,...,l_r)$. If
$l_1+...+l_r=n+1$ then $X$ is a Calabi-Yau manifold and its quantum cohomology
is described by the mirror conjecture. In this and the next sections we study
respectively the cases $l_1+...+l_r < n$ and $l_1+...+l_r=n$
when the $1$-st Chern class of $X$ is still positive. In the case $l_1+...+l_r
> n+1$ (which from the point of view of enumerative geometry can be considered
as ``less interesting'' for rational curves generically occur only in
finitely many degrees) the ``mirror symmetry'' problem of hypergeometric
interpretation of quantum cohomology differential equations remains open.

Let $E_d$ be the Euler class of the vector bundle over the moduli space
$Y_{2,d}$ of genus $0$ degree $d$ stable maps $\phi : (C, x_0,x_1) \to \CC P^n$
with two marked points, with the fiber $H^0(C, \phi ^* H^{l_1}\oplus ...
\oplus \phi^* H^{l_r})$ where $H^l$ is the $l$-th tensor power
of the hyperplane line bundle over $\CC P^n$.

Consider the class

$$S_d(\hbar ):=\frac{1}{\hbar + c_1^{(0)}}\ E_d  \ \in \ H^*(Y_{2,d})$$
where $c_1^{(0)}$ is the $1$-st Chern class of
the ``universal tangent line at the marked point $x_0$ '',
and $e_0, e_1$ are the evaluation maps. Due to the factor $E_d$ this
class represents the push forward along $X_{2,d} \to Y_{2,d}$ of the class
$1 /(\hbar + c_1^{(0)}) \ \in H^*(X_{2,d})$
(by the very construction of $X_{2,d}$ in Section $4$).

In the cohomology algebra $\CC [P]/(P^{n+1})$ of $\CC P^n$, consider the class

\[ S(t,\hbar):= e^{Pt/\hbar }\sum _{d=0}^{\infty } e^{dt}
(e_0)_* (S_d(\hbar )) \]
where $(e_0)_*$ represents the push-forward along the evaluation map
(and for $d=0$, when $Y_{2,d}$ is not defined, we take
$Euler \ (\oplus _j H^{\otimes l_j}) $ on the role of $(e_0)_* S_0$).

Considered as a function of $t$, $S$ is a curve in $H^*(\CC P^n)$ whose
components are solutions of the differential equation we are concerned
about. Indeed, according to Section $6$ a similar sum represents
the solutions of the quantum cohomology differential equation for $X$, and
$S$ is just the push-forward of that sum from $H^*(X)$ to $H^*(Y)$.
(Strictly speaking $S$ carries information only about correlators between
those classes which come from the ambient projective space; also if
$X$ is a surface $\rk H_2(X)$ can be greater than $1$ and $S$ mixes information
about the curves of different degrees in $X$ when they have the same degree
in $Y$.)

\bigskip

{\bf Theorem 9.1.} {\em Suppose that $l_1+...+l_r < n$. Then}

\[ S = e^{Pt/\hbar } \sum _{d=0}^{\infty } e^{dt}
\frac{ \Pi _{m=0}^{dl_1} (l_1P + m\hbar) ... \Pi _{m=0}^{dl_r} (l_rP+m\hbar)}
{\Pi _{m=1}^d (P+m\hbar )^{n+1} } .\]

The formula coincides with those in \cite{HG1, HG} (found by analysis
of toric compactifications of spaces of maps $\CC P^1 \to \CC P^n$) for
solutions of differential equations in $S^1$-equivariant cohomology of the
loop space.

\bigskip

{\bf Corollary 9.2.} (see \cite{HG1, HG})
 {\em The components $s:=\lan P^i, S\ran , i=0,...,n-r, $
of $S$ form a basis of solutions to the linear differential equation}

\[ (\hbar \frac{d}{dt})^{n+1-r}\ s \ = \
e^t \Pi _{j=1}^r \ l_j\ \Pi _{m=1}^{l_j-1}\
\hbar (l_j\frac{d}{dt} + 1) \ s \ .\]

This implies (combine \cite{HG} with \cite{BVS}) that the solutions
have an integral representation of the form

\[ \int _{\g ^{n-r}\subset X'_t} \ e^{(u_0+...+u_n)/\hbar } \ \
\frac{du_0\w ... \w du_n}{dF_0\w dF_1\w ... \w dF_r} \]
where
\[ F_0=u_0...u_n,\ F_1=u_1+...u_{l_1},\ F_2=u_{l_1+1}+...+u_{l_1+l_2},
... , \ F_r=u_{l_1+...+l_{r-1}+1}+...+u_{l_1+...+u_{l_r}} \]
and the ``mirror manifolds'' $X'_t$ are described by the equations
\[ X'_t=\{ (u_o,...,u_n)\ |\ F_0(u)=e^t,\ F_1(u)=1, ... ,\ F_r(u)=1 \} .\]
%
%
This proves for $X$ the mirror conjecture in the form suggested in
\cite{HG}.

\bigskip

{\bf Corollary 9.3.} {\em If $\dim _{\CC } X \neq 2$ the cohomology class
$p$ of hyperplane section
satisfies in the quantum cohomology of $X$ the relation}

\[ p^{n+1-r}=l_1^{l_1}...l_r^{l_r} q p^{l_1+...+l_r-r} .\]

When $X$ is a surface the same relation holds true in the quotient of the
quantum cohomology algebra which takes in account only degrees of curves in
the ambient $\CC P^n$ (we leave to figure out a precise description of this
quotient to the reader; quadrics $\CC P^1\times \CC P^1$ in $\CC P^3$ provide
a good example:
$(p_1+p_2)^3 = 4q(p_1+p_2) \ \mod p_1^2=q=p_2^2 $.)

This corollary is consistent with the result of A. Beauville \cite{Bea}
describing quantum cohomology of complete intersections with
$\sum l_j \leq n+1 -\sum (l_j-1) $ and with results of M. Jinzenji \cite{J} on
quantum cohomology of projective hypersurfaces ($r=1$) with $l_1<n$.

\bigskip

{\bf Corollary 9.4.} {\em The number of degree $d$ holomorphic maps
$\CC P^1 \to X^{n-r} \subset \CC P^n$,
which send $0$ and $\infty $ to two given
cycles and send $n+1-r$ given points in $\CC P^1$  to
$n+1-r$ given generic hyperplane sections, is equal to
$l_1^{l_1}...l_r^{l_r}$ times
the number of degree $d-1$ maps which send $0$ and $\infty $ to the same cycles
and $l_1+...+l_r-r$ given points --- to $l_1+...+l_r-r$ given
hyperplane sections.}

This is the enumerative meaning of Corollary $9.3$; of course in this
formulation numerous general position reservations are assumed.

{\em Control examples.} $1.$ $l_1=...=l_r=1$: The above formulas for quantum
cohomology
and for solutions of the differential equations in the case of a hyperplane
section give rise to the same formulas with $n:=n-1$.

$2$. $n=5, r=1, l=2$: $X$ is the Plucker embedding of the grassmanian
$Gr_{4,2}$. Its quantum cohomology algebra is described by the relations
$c_1^3=2c_1 c_2,\ c_2^2-c_2c_1^2+q=0$ between the Chern classes of the
tautological plane bundle. For the $1$-st Chern class $p=-c_1$ of the
determinant line bundle we deduce the relation $p^5=4pq$ prescribed by
Corollary $9.3$.

\bigskip

We will deduce Theorem $9.1$ from its equivariant generalization.
Consider the space $\CC ^{n+1}$ provided with the standard action of the
$(n+1)$-dimensional torus $T$. The equivariant cohomology algebra
of $\CC ^{n+1}$ coincides with the algebra of characteristic classes
$H^*(BT^{n+1}) = \CC [\l _0,...,\l _n]$. The equivariant cohomology algebra
of the projective space $(\CC ^{n+1} - 0)/\CC ^{\times} $ in these notations
is identified with $\CC [p,\l ]/ ((p-\l_0)...(p-\l_n))$ and the push-forward
$H_T^*(\CC P^n)\to H_T^*(pt)$ is given by the residue formula

\[ f(p,\l ) \mapsto \frac{1}{2\pi i} \oint
\ \frac{f(p,\l )dp}{(p-\l_0)...(p-\l_n)} \ .\]

Here $-p$ can be considered as the equivariant $1$-st Chern class of the
Hopf line bundle provided with the natural lifting of the torus action.
We will use $\phi _i:=\Pi _{j\neq i} (p-\l _j), \ i=0,...,n, $ as a basis in
$H_T^*(\CC P^n)$.

Consider the $T$-equivariant vector bundle $\oplus _{j=1}^r H^{\otimes l_j}$
and provide it with the fiberwise action of the additional $r$-dimensional
torus $T'$. The equivariant Euler class of this bundle is equal to
$ (l_1p-\l '_1)...(l_rp-\l '_r)$ where $\CC [\l ']=H^*(BT')$.

Introduce the equivariant counterpart $S'$ of the class $S$ in the
$T\times T'$-equivariant cohomology of $\CC P^n$. This means that we use
the equivariant class $p$ instead of $P$ and replace the Euler classes
$E_d$ and $c_1^{(0)}$ by their equivariant partners.

\bigskip

{\bf Theorem 9.5.} {Let $l_1+...+l_r < n$. Then}

\[ S'= e^{pt/\hbar} \sum _{d=0}^{\infty }
e^{dt} \frac{ \Pi _0^{dl_1} (l_1p-\l '_1 +m\hbar) ...\Pi _0^{dl_r}
(l_rp-\l '_r +m\hbar)}{\Pi _1^d (p-\l _0+m\hbar) ... \Pi _1^d (p-\l _n+m\hbar)}
\ .\]

Theorem $9.1$ follows from Theorem $9.5$ by putting $\l =0, \l '=0$ which
corresponds to passing from equivariant to non-equivariant cohomology.

The vector-function $S'$ satisfies the differential equation

\[ \Pi _{i=0}^r (\hbar \frac{d}{dt} -\l _i) \ S' =
e^t \Pi _{m=1}^{l_1} (l_1\hbar \frac{d}{dt} - \l '_1 + m\hbar ) ...
    \Pi _{m=1}^{l_r} (l_r\hbar \frac{d}{dt} - \l '_r + m\hbar ) \ S' .\]

%
%

\bigskip

We intend to prove Theorem $9.5$ by means of localization of $S'$ to the
fixed point set of the torus $T$ action on the moduli spaces $Y_{2,d}$.
As it is shown in \cite{Kn}, all correlators of the equivariant theory on
$\CC P^n$ are computable at least in principle, and in practice the
computation reduces to a recursive procedure which can be understood as a
summation over trees and can be also formulated as a non-linear fixed point
(or critical point) problem. We will see below that in the case of correlators
$\lan \phi _i, S'\ran $ certain reasons of a somewhat geometrical character
cause numerous cancellations between trees so that the recursive procedure
reduces to a ``summation over chains'' and respectively to a {\em linear}
recurrence equation. The formula of Theorem $9.5$ is simply the solution to
this equation.

\bigskip

In the proof of Theorem $9.5$ below we write down all formulas for for $r=1$
(it serves the case when $X$ is a hypersurface in $\CC P^n$ of degree $l < n$).
The proof for $r > 1$ differs only by longer product formulas.

\bigskip

Let us abbreviate $c_1^{(0)}$ as $c$,
denote $E'_d$ the equivariant Euler class of the vector bundle over $Y_{2,d}$
whose fiber over the point $\psi : (C, x_0, x_1) \to Y=\CC P^n$ consists of
holomorphic sections of the bundle $ \psi ^* (H^l)$
{\em vanishing at} $x_0$, and introduce the following equivariant correlator:

\[ Z_i :=\sum _{d=0}^{\infty } q^d
\int _{Y_{2,d}} e_0^* (\phi _i) \frac{1}{\hbar + c} E'_d .\]

We have

\[ \lan \phi _i, S' \ran = e^{\l _i t/\hbar } (l\l _i - \l ') (Z_i | _{q=e^t})
\]

\medskip

{\bf Proposition 9.6.}
\[ Z_i = 1 + \sum _{d>0} (\frac{q}{\hbar ^{n+1-l}})^d
\int _{Y_{2,d}} \frac{(-c)^{(n+1-l)d-1}}{1+c/\hbar} \ E'_d \
e_0^*(\phi _i) \ . \]

{\em Proof.} We have just dropped first several terms in the geometrical
series $1/(\hbar + c)$ since their degree added with the degrees of other
factors in the integral over $Y_{2,d}$ is still less than the dimension
of $Y_{2,d}$. It is important here that all the equivariant classes involved
including $\phi _i$ are defined in the equivariant cohomology over
$\CC [\l ,\l ']$ without any localization.

\medskip

It is a half of the geometrical argument mentioned above. The other half
comes from the description of the fixed point set in $Y_{2,d}$ given
in \cite{K}.

Consider a fixed point of the torus $T$ action on $Y_{2,d}$. It is represented
by a holomorphic map of a possibly reducible curve with complicated
combinatorial structure and with two marked points on some components. Each
component carrying $3$ or more special points is mapped to one of the $n+1$
fixed points of $T$ on $\CC P^n$, and the other components are mapped
(with some multiplicity) onto the lines joining the fixed points and
connect the point-mapped components in a tree-like manner.

In the Borel localization formula for $\int e_0^*(\phi _i) ...$ the fixed
point will have zero contribution unless the marked point $x_0$ is mapped
to the $i$-th fixed point in $\CC P^n$ (since $\phi _i$ has zero localizations
at all other fixed points.

Consider a fixed point curve $C$ whose marked point $x_0$ is indeed mapped to
the $i$-th fixed point in $\CC P^n$. There are two options

(i) the marked point $x_0$ is situated on an irreducible component of $C$
mapped with some degree $d'$ onto the line joining the $i$-th
fixed point with the $j$-th fixed point in $\CC P^n$ with $i\neq j$;

(ii) the marked point $x_0$ is situated on a component of $C$ mapped to
the $i$-th fixed point and carrying two or more other special points.

Consider first the option (ii) and the contribution of such a connected
component of the fixed point set in $Y_{2,d}$ to the Borel localization
formula for $\int c^{(n+1-l)d-1} ...$. The connected component itself
is the (product of the) Deligne-Mumford configuration space of, say, $s+1$
special points: the marked point $x_0$, may be the marked point $x_1$, and
respectively $s-1$ or $s$ endpoints of other components of $C$ mapped
onto the lines outgoing the $i$-th fixed point in $\CC P^n$.

{\bf Lemma 9.7.} {\em The type (ii) fixed point component
in $Y_{2,d}$ has zero contribution to the Borel localization formula for}
$\int _{Y_{2,d}} c^{(n+1-l)d-1} ... $

{\em Proof.} Restriction of the class $c$ from $Y_{2,d}$ to the type (ii)
fixed point component coincides with the $1$-st Chern class of the line
bundle on the Deligne-Mumford factor $\calm _{0, s+1}$ of the component
defined as ``the universal tangent line as the marked point $x_0$''
and is thus
nilpotent in the cohomology of the component. Since the number of straight
lines in a curve of degree $d$ does not exceed $d$ we find that the dimension
$s-2$ of the factor $\calm _{0,s+1}$ is less than $d$ which in its turn does
not exceed $(n+1-l)d-1$ for $d>0$ (because we assumed that $n+1-l\geq 2$).

\medskip

Consider now the option (i). The irreducible component $C'$ of the curve
$C=C'\cup C''$
carrying the marked point $x_0$ is mapped with the multiplicity $d' \leq d$
onto the line joining $i$-th fixed point in $\CC P^n$ with the $j$-th one
{\em while the remaining part $C''\to \CC P^n$ of the map represents a fixed
point in $Y_{2,d-d'}$}. Moreover, the normal space to the fixed point
component at the type (i) point (the equivariant Euler class of the normal
bundle occurs in the denominator of the Borel localization formula) is the
sum of (a) such a space $N''$ for $C'' \to \CC P^n$,  (b) the space $N'$
of holomorphic vector fields along the map $C'\to \CC P^n$ vanishing at
the fixed point $j$ factorized by infinitesimal reparametrizations of $C'$,
(c) the tensor product $L$ of the tangent lines to
$C'$ and $C''$ at their intersection point. Since the space $V$ of holomorphic
sections of $H^l$ restricted to $C$ (and vanishing at $x_0$) admits a similar
decomposition $V'\oplus V''$, we arrive to the following linear recursion
relation for $Z_i$.

{\bf Proposition 9.8.} {\em Put $z_i(Q,\hbar ):=
Z_i(\hbar ^{(n+1-l)}Q, \hbar )$.  Then}

\[ z_i(Q, \hbar) = 1 + \sum _{j\neq i} \sum _{d'>0} Q^{d'} Coeff \ _i^j (d')
\ z_j(Q, (\l _j-\l _i)/d' ) \]
{\em where}
\[ Coeff \ _i^j(d')=
\frac{[(\l_j -\l_i)/d]^{(n+1-l)d'-1}}{1+(\l _i-\l _j)/d' \hbar}
\frac{Euler \ (V')}{Euler \ (N')} \ \phi _i |_{p=\l _i} \ . \]

{\em Proof.} Here $(\l _i-\l_j)/d' $ is the localization of $c$, and the key
point is that the equivariant Chern class of the line bundle $L$ over
$Y_{2,d-d'}$ is what we would denote $\hbar + c$ {\em for the moduli space}
$Y_{2,d-d'}$ but with $\hbar = (\l _j -\l _i)/d'$. This is how the recursion
for the correlators $z_i$ becomes possible. The rest is straightforward.

{\em Remark.} Our reduction to the linear
recursion relation can be interpreted in the following more geometrical way:
contributions
of all non-isolated fixed points cancel out with some explicit part of the
contribution from {\em isolated} fixed points; the latter are represented by
chains of multiple covers of straight lines connecting the two marked points.

\bigskip

Let us write down explicitly the factor $Coeff \ _i^j(d)$
from Proposition $9.8$
(compare with \cite{Kn}). $Coeff \ _i^j (d) =$

\[ \frac{
\Pi _{m=1}^{ld} (l\l _i -\l' + m(\l _j-\l_i)/d) [(\l_j-\l_i)/d]^{(n+1-l)d-1}
}{
d (1+(\l_i -\l_j)/\hbar d)
\Pi _{\a =0}^n \ _{m=1}^d \ _{(\a ,m)\neq (j, d)} (\l_i-\l_{\a}+m(\l_j-\l_i)/d)
} = \]
(here the product in the numerator is $Euler (V')$, the
denominator --- it is essentially $Euler (N')$ where however the cancellation
with $\phi _i |_{p=\l _i}$ is taken care of ---  has been computed using the
exact sequence $0 \to \CC \to \CC ^{n+1} \otimes H \to T_Y \to 0$ of vector
bundles over $Y= \CC P^n$, and the ``hard-to-explain'' extra-factor $d$
is due to the orbifold structure of the moduli spaces (the $d$-multiple map
of $C'$ onto the $(ij)$-line in $\CC P^n$ has a discrete symmetry of order
$d$)

\[ = \frac{1}{[(\l_i-\l_j)/\hbar + d]}
\frac{\Pi _{m=1}^{ld} (\frac{(\l_i-\l')d}{\l_j-\l_i} +m)}
{\Pi _{\a =0}^n\ _{m=1}^d \ _{(\a ,m)\neq (j,d)}
(\frac{(\l_i-\l_{\a })d}{\l_j-\l_i} + m) } \ . \]

Now it is easy to check

{\bf Proposition 9.9.} {\em The correlators $z_i(Q, 1/\o )$ are power series
$\sum C_i(d) Q^d $ in $Q$ with coefficients $C_i(d)$ which are {\em reduced}
rational functions of $\o $ with poles of the order $\leq 1$ at
$\o = d'/(\l_j-\l_i)$ with $d'=1,...,d$. The correlators $z_i$ are uniquely
determined by these properties, the recursion relations of Proposition $9.8$
and the initial condition $C_i(0)=1$.}

\bigskip

The proof of Theorem $9.5$ is completed by the following

{\bf Proposition 9.10.} {\em The series}

\[ z_i=\sum _{d=0}^{\infty } Q^d
\frac{\Pi _{m=1}^{ld} ((l\l_i -\l')\o +m)}
{d! \Pi _{\a \neq i} \Pi _{m=1}^d ((\l_i -\l_{\a })\o +m)} \]

{\em satisfy all the conditions of Proposition $9.10$.}

{\em Proof.} The recursion relation is deduced by the decomposition of
the rational functions of $\o $ into the sum of simple fractions (or,
equivalently, from the Lagrange interpolation formula for each numerator
through its values at the roots of the corresponding denominator).

\section{Complete intersections with $l_1+...+l_r=n$}
\label{sec10}

Let $X\subset Y=\CC P^n$ be a non-singular complete intersection given by
equations of degrees $(l_1,...,l_r)$ with $l_1+...+l_r=n$. There are only
two points where our proof of Theorem $9.1$ would fail for such $X$. One of
them is the Lagrange interpolation formula in the proof of Proposition $9.10$.
Namely, the rational functions of $\o $ there are not reduced ---
the degree $dl$ of the numerator {\em is equal} to the degree $dn$ of the
corresponding denominator. The other one is Lemma $9.7$. Namely, we have the
following lemma instead.

{\bf Lemma 10.1.}
{\em The type (ii) fixed point component in $Y_{2,d}$ makes zero
contribution via Borel localization formulas to $\int _{Y_{2,d}} c^{d-1}...$
{\em unless} it consists of maps $(C'\cup C'', x_0,x_1)\to Y$ where $C'$ is
mapped to a fixed point in $\CC P^n$ and carries both marked points, and
$C''$ is a disjoint union of $d$ irreducible components (intersecting $C'$ at
$d$ special points) mapped (each with multiplicity $1$) onto straight lines
outgoing the fixed point. All type (ii) components make zero contribution
to $\int _{Y_{2,d}} c^d ...$.}

Let us modify the results of Section $9$ accordingly. As we will see,
the LHS in Theorem $9.5$ is now only {\em proportional} to the RHS, and we will
compute the proportionality coefficient (a series in $q$) directly.

{\bf Proposition 10.2.} {\em Put $z_i(Q,\hbar ):= Z_i(\hbar Q,\hbar )$. Then}

\[ z_i(Q,\hbar) = 1 + \sum _{d>0} Q^d Coeff\ _i(d) +
 \sum _{j\neq i} \sum _{d'>0} Q^{d'}
Coeff \ _i^j(d')\ z_j(Q,(\l_j-\l_i)/d') \]
{\em where $Coeff \ _i(d)$ is equal to the contribution of type (ii) fixed
point components to \newline
$\int _{Y_{2,d}} (-c)^{d-1} E_d' \ e_0^*(\phi _i) $, and}
\[ Coeff \ _i^j(d) =  \frac{1}{[(\l_i-\l_j)/\hbar +d]}\
\frac{\Pi _{a=1}^r \Pi _{m=1}^{dl_a} (\frac{(l_a\l_i-\l_a')d}{\l_j-\l_i}+m)}
{\Pi _{\a =0}^n \ _{m=1}^d \ _{(\a ,m)\neq (j,d)}
(\frac{(\l_i-\l_{\a })d}{\l_j-\l_i} +m)} \ .\]

{\bf Corollary 10.3.} {\em The correlators $z_i(Q,1/\o )$ are power series
$\sum _d C_i(d) Q^d$ with coefficients
\[ C_i(d)=P_d(\o ,\l ,\l ')/\Pi _{\a }\Pi _{m=1}^d ((\l_i-\l_{\a })\o +m) \]
where $P_d = P_d^0 \o ^{nd} + ... $ is a polynomial in $\o $ of degree $nd$.
The correlators $z_i$ are uniquely determined by these properties, the
recursion relations of Proposition $10.2$ and the initial conditions
\[ \sum _d Coeff \ _i(d) Q^d = \sum _d Q^d \frac{P_d^0}
{d!\Pi _{\a \neq _i} (\l _i-\l_{a })^d}  \ .\]}

{\bf Proposition 10.4.} {\em The series
\[ z_i' =\sum _{d=0}^{\infty } Q^d
\frac{\Pi _{a=1}^r \Pi _{m=1}^{l_a d} ((l_a\l_i-\l_a')\o +m)}
{d!\Pi _{\a \neq i}\Pi _{m=1}^d ((\l_i-\l_{\a })\o +m)} \ \]
satisfy the requirements of Corollary $10.3$ with the initial condition
\[ \sum _d Q^d \frac{\Pi _{a=1}^r (l_a\l_i-\l_a')^{l_a d} }
{ d!\Pi _{\a \neq i}
(\l_i-\l_{\a })^d} \ = \
 \exp \{ Q \frac{\Pi _a (l_a\l_i -\l_a')^{l_a}}{\Pi _{\a \neq i}
(\l_i-\l_{\a })} \} . \]}

\bigskip

Now let us compute $Coeff \ _i(d)$ using the description of type (ii) fixed
point components given in Lemma $10.1$.

{\bf Proposition 10.5.} {\em Contribution of the type (ii) fixed point
components to \newline
$\sum _d Q^d \int _{Y_{2,d}} (-c)^{d-1} E'_d \ \phi _i $ is}
\[ \exp \{ Q \frac{\Pi _a (l_a\l_i-\l_a')^{l_a}}
{\Pi _{\a \neq i} (\l_i-\l_{\a })} \} \ \exp \{ - Q \ l_1!...l_r!  \} \ .\]

{\em Proof.} Each fixed point component described in Lemma $10.1$ is isomorphic
to the Deligne - Mumford configuration space $\calm _{0,d+2}$. Our computation
is based on the following known formula (see for instance \cite{Kn} ) for
correlators between Chern classes of universal tangent lines at the marked
points:
\[ \int _{\calm _{0, k}} \frac{1}{(w_1+c_1^{(1)}) ... (w_k+c_1^{(k)})} =
\frac{(1/w_1+...1/w_k)^{k-3}}{w_1 ... w_k} \ .\]

Consider the type (ii) fixed point component specified by
the following combinatorial structure of stable maps: $d$ degree $1$
irreducible components join the $i$-th fixed point with the fixed points
with indices $j_1,...,j_d$. Using the above formula and describing explicitly
the normal bundle to this component in $Y_{2,d}$ and localization of the Euler
class $E_d'$ we arrive to the following expression for the contribution of this
component to $\int _{Y_{2,d}} (-c)^{d-1} \phi _i E_d'$:
\[  \Pi _{s=1}^d
\frac{\Pi _{a=1}^r \Pi _{m=1}^{l_a}(l_a\l_i-\l_a' +m(\l_{j_s}-\l_i))}
{(\l_i-\l_{j_s}) \Pi _{\a \neq j_s,i} (\l_{j_s}-\l_{\a })} \ . \]

Summation over all type (ii) components in all $Y_{2,d}$ with weights $Q^d$
gives rise to

\[ \exp \{ - Q \sum _{j\neq i}
\frac{\Pi _a \Pi _{m=1}^{l_a} (l_a\l_i-\l_a' +m(p-\l_i))}
{\Pi _{\a \neq j} (p-\l_{\a})} \ |_{p=\l_j}  \} \ .\]
The exponent can be understood as a sum of residues at $p\neq \l_i, \infty $
and is thus opposite to the sum
\[ l_1!...l_r! -
\frac{\Pi _a (l_a\l_i-\l_a')^{l_a}}{\Pi _{\a \neq i} (\l_i-\l_{\a })}
 \]
of residues at $\infty $ and $\l_i $.

{\bf Corollary 10.6.}  $z_i(Q,1/\o )= z_i'(Q,\o )\ \exp (- l_1! ... l_r! Q)$.

{\em Proof.} Multiplication by a function of $Q$ does not destroy the
recursion relation of Proposition $10.2$ but changes the initial condition.

\bigskip

We have proved the following

{\bf Theorem 10.7.} {\em Suppose $l_1+...+l_r=n$. Then}

\[ S' = e^{ (pt - l_1!...l_r! e^t)/\hbar } \
\sum _{d=0}^{\infty } \ e^{dt} \
\frac{\Pi _0^{dl_1} (l_1p-\l_1'+m\hbar) ... \Pi _0^{dl_r} (l_rp-\l_r'+m\hbar)}
{\Pi _1^d (p-\l_0 +m\hbar ) ... \Pi _i^d (p-\l_n +m\hbar )} \ .\]

\[ S = S' |_{\l =0, \l' =0} =
e^{ (Pt - l_1! ... l_r! e^t)/\hbar }
\frac{ \Pi _{j=1}^r \Pi _{m=0}^{dl_j} (l_j P+m\hbar)}
{\Pi _{m=1}^d (P+m\hbar)^{n+1}} \ \ (\text{mod} \ P^{n+1}) \ .\]

{\bf Corollary 10.8.} {\em Let $D=\hbar d/dt + l_1!...l_r! e^t$. Then}
\[ D^{n+1-r} S = l_1...l_r e^t \Pi _{j=1}^r
(l_j D+\hbar)...(l_j D +(l_j-1)\hbar ) \ S .\]

{\bf Corollary 10.9.} {\em In the quantum cohomology algebra of $X$ the class
$p$ of hyperplane sections satisfies the following relation
(with the same reservation in the case $\dim X \leq 2$ as in Corollary $9.3$):}
\[ (p + l_1!...l_r! q)^{n+1-r}=
l_1^{l_1} ... l_r^{l_r} q (p+l_1!...l_r!q)^{n-r} \ .\]

{\em Control examples.} $1$.
$X=pt$ in $\CC P^1$ ($n=1, r=1, l=1$). The above relation
takes on $p+q=q$, or $p=0$. Since $P^2=0$, we also find from Theorem 10.7 that
$ S=P \exp (-e^t) \sum _d e^{dt}/d! = P $, or $\lan 1, S\ran =1 $ as it should
be for the solution of the differential equation $\hbar d/dt \ s = 0$ that
arises from quantum cohomology of the point.

$2$. $X=\CC P^1$ embedded as a quadric into $\CC P^2$ ($n=2, r=1, l=2$).
We get $(p+2q)^2=4q(p+2q)$, or $p^2=4q^2$. Taking into account that $p$ is
twice the generator in $H^2(\CC P^1)$ and the line in $\CC P^1$ has the degree
$2$ in $\CC P^2$ we conclude that this is the correct relation in the quantum
cohomology of $\CC P^1$. This example was the most confusing for the author:
predictions of the loop space analysis \cite{HG1} appeared totally unreliable
because they gave a wrong answer for the quadric in $\CC P^2$. As we see now,
the loop space approach gives correct results if $l_1+...+l_r<n$ and require
``minor'' modification (by the factor $\exp (-l_1!...l_r! q /\hbar )$ )
in the boundary
cases $l_1+...l_r=n$; the quadric on the plane happens to be one of such cases.


$3$. $n=3, r=1, l=3$. We have $(p+6q)^3 = 27q (p+6q)^2 $, or $p^3=9qp^2+6^3q^2p
+ 27\cdot 28 q^3$. In particular, $\lan p * p, p\ran = 9 q \lan p, p \ran +
6^3q^2\lan p, 1\ran + 27\cdot 28 q^3 \lan 1, 1\ran = 27 q +0+0$ which indicates
that there should exist $27$ discrete lines on a generic cubical surface in
$\CC P^3$.

\section{Calabi-Yau projective complete intersections}
\label{sec11}

Let $L_d(Y)$ denote, as in Section $6$, the moduli space of stable maps
$\psi: \CC P^1 \to \CC P^n \times \CC P^1$ of bidegree $(d,1)$ with $2$ marked
points mapped to $\CC P^n \times \{ 0\} $ and $\CC P^n \times \{ \infty \} $
respectively. Let ${\cal E} _d$ denote the equivariant Euler class of the
vector bundle over $L_d(Y)$ with the fiber
$H^0(\CC P^1, \psi ^*(E))$ where $E$ is the bundle on $\CC P^n \times \CC P^1$
induced from our ample bundle $\oplus _a H^{l_a}$ by the projection to the
first factor.

Consider the equivariant correlator

\[ \Phi = \int _Y \ Euler^{-1}(E) \ S'(t,\hbar ) \ S'(\t ,-\hbar) \ = \]
\[ = \sum _{d,d'} e^{dt} e^{d'\t} \sum _i \frac{\Pi _a (l_a \l_i -\l'_a)}
{\Pi _{j\neq i} (\l_i-\l_j)}
\int _{Y_{2,d}} E'_d \frac{e^{pt/\hbar } e_0^*(\phi _i)}{\hbar + c} \
\int _{Y_{2,d'}} E'_{d'} \frac{e^{-p\t /\hbar} e_0^*(\phi _i)}{-\hbar +c} .\]

In the case $l_1+...+l_r < n$ it is easy to check using the explicit formula
for $S'$ from Theorem $9.5$ that
\[ \Phi = \frac{1}{2\pi i} \oint e^{p(t-\t )/\hbar}  [ \sum _d e^{d\t }
\frac{\Pi _a \Pi _{m=0}^{l_a d} (l_a p-\l'_a -m\hbar)}
{\Pi _{j=0}^n \Pi _{m=0}^d (p-\l_j - m\hbar)}] \ dp.  \]

This is an equivariant version of a formula found in \cite{HG1} in the context
of loop spaces and toric compactifications of spaces of rational maps. Namely,
consider the projective space $L'_d$ of $(n+1)$-tuples of polynomials in
one variable of degree $\leq d$ each, up to a scalar factor (notice that
$L'_d$ has the same dimension $d(n+1)+n$ as $L_d$). It inherits the
component-wise action of the torus $T^{n+1}$ and the action of
$S^1$ by the rotation
of the variable (``rotation of loops''). Integration over the equivariant
fundamental cycle in $L'_d$ is given by the residue formula

\[ f(p,\l , \hbar) \mapsto \frac{1}{2\pi i}\oint \frac{ f dp }
{\Pi _{j=0}^n \Pi _{m=0}^d (p -\l_j - m\hbar ) } .\]

Consider the equivariant vector bundle over $L'_d$ such that substitution
of the $(n+1)$ polynomials into $r$ (invariant) homogeneous equations
in $\CC P^n$ of degrees $l_1, ..., l_r$ produces a section of this bundle.
The equivariant Euler class of the bundle is
\[ {\cal E} '_d = \Pi _{a=1}^r \Pi _{m=0}^{l_a d} (l_a p -\l'_a - m\hbar ) .\]

The formula for $\Phi $ indicates that there should exist a close relation
between the spaces $L_d$ and $L'_d$. This relation is described in the
following lemma whose proof will be given in the end of this Section.

{\bf The Main Lemma.} {\em There exists a natural
$S^1\times T^{n+1}$-equivariant map
$\m : L_d\to L'_d$. Denote $-p$ the equivariant $1$-st Chern class of the Hopf
bundle over $L'_d$ induced by $\m $ to $L_d$. Then}

\[ \Phi (t,\t) =
\sum _d e^{d\t } \int _{L_d} e^{p(t-\t )/\hbar} {\cal E} _d .\]

Define $\Phi ' (q, z, \hbar ) := \Phi | t=\t +z\hbar, q=e^{\t}$. (By the way
the limit of the series $\Phi '$ at $\hbar =0$ has the topological meaning
of what is called in \cite{GK} the {\em generating volume function}, and the
meaning of this limit procedure in terms of differential equations satisfied
by $\Phi $ is the {\em adiabatic approximation}.)

{\bf Corollary 11.1.} $\Phi '(q,z):= \sum_d q^d \int _{L_d} e^{pz} {\cal E} _d
 =$

\[ = \frac{1}{2\pi i} \oint e^{pz}  \sum _d \frac{ E_d(p,\l, \l',\hbar )}
{\Pi _{j=0}^n \Pi _{m=0}^d (p-\l_j -m\hbar)} dp \]
{\em where $E_d = \m _* ({\cal E} _d)$ is a {\em polynomial} (of degree
$<(n+1)d$) of all its variables.}

{\em Proof.} The integrals $E^{(k)} = \int _{L_d} p^k {\cal E} _d, \ k=0,...,
\dim L'_d$, which determine the push-forward $\m _* ({\cal E})$
are polynomials in $(\l, \l', \hbar )$. The matrix
$\int _{L'_d} p^{i+\dim L'_d - j}$ is triangular with all eigenvalues equal
to $1$. This means that there exists a unique polynomial in $p$ with
coefficients {\em polynomial in} $(\l, \l', \hbar)$ which represents
the push-forward with any given polynomials $E^{(k)}(\l ,\l', \hbar)$.

The last argument also proves

{\bf Proposition 11.2.}{ \em Suppose that a series
\[ s= \sum _d q^d \frac{ P_d(p,\l,\l',\hbar)}
{ \Pi _j \Pi _{m =0}^d (p-\l_j -m\hbar)} \]
with coefficients $P_d$ which are polynomials of $p$ of degree $\leq \dim L_d$
has the property that for every $k=0,1,2,...$ the $q$-series $\oint s p^k dp $
has polynomial coefficients in $(\l, \l', \hbar)$. Then the coefficients of
all $P_d$ are polynomials of $(\l,\l',\hbar)$, and vice versa.}

\bigskip

The coefficient $E_d(p,\l,\l',\hbar)$ in the series $\Phi '$ has the total
degree $(l_1+...+l_r)d+r$ according to the dimension of the vector bundle
whose Euler class it represents. Consider the following operations with the
series $\Phi $:

(i) multiplication by a series of $e^t$ and / or $e^{\t }$;

(ii) simultaneous change of variables $t\mapsto t+f(e^t),
\t \mapsto \t +f(e^{\t} )$.

(iii) multiplication by $\exp [C (f(e^t)-f(e^{\t }))/\hbar ]$
(here the factor $C$ should be a linear function of $(\l, \l')$ in order
to obey homogeneity).

{\bf Proposition 11.3.} {\em  The property of the series $\Psi $ to generate
polynomial coefficients $E_d(p,\l,\l' \hbar)$ is invariant with respect
to the operations (i),(ii),(iii).}

{\em Proof.} The polynomiality property of coefficients in $\Phi '$ is
equivalent, due to Proposition $11.2$, to the fact that for all $k$ the
$q$-series $(\p /\p z)^k |_{z=0} \Phi ' $ has polynomial coefficients.

Multiplication by a series of $q$ does not change this property, which proves
the invariance with respect to multiplication by functions of $e^{\t }$.

The roles of $t$ and $\t $ can be interchanged by the substitutions
$p \mapsto p+\hbar d, \hbar \mapsto -\hbar $ in each summand of $\Phi $.
This proves the invariance with respect to multiplication by functions
of $e^t$.

The operation (ii) transforms $\Phi '$ to
\[ \frac{1}{2\pi i} \sum _d q^d e^{d f(q)} \oint  \
\exp \{ p\frac{z\hbar + f(qe^{z\hbar})-f(q)}{\hbar } \} \
\frac{E_d(p)}{\Pi _j \Pi _m (p-\l_j -m\hbar)} \ dp .\]

Since the exponent is in fact divisible by $\hbar $, the derivatives in $z$
at $z=0$ still have polynomial coefficients. This proves the invariance
with respect to (ii). The case of the operation (iii) is analogous.

\bigskip

  We are going to use the above polynomiality and invariance
properties of the correlator $\Phi $ in order to describe quantum cohomology
of {\em Calabi-Yau} complete intersections in $\CC P^n$ (in which case
$l_1+...+l_r = n+1$). We will use this polynomiality in conjunction
with recursion relations based on
the fixed point analysis of Sections $9, 10$. The result can be roughly
formulated in the following way: the hypergeometric functions of Theorem $9.5$
in the case $l_1+...+l_r =n+1$ can be transformed to the correlators $S'$ by
the operations (i),(ii),(iii). Notice that in the Calabi -- Yau case
all our formulas are homogeneous with the grading $\deg q = 0,
\deg p =\deg \hbar =\deg \l =\deg \l' =1, \deg z =-1$. In particular the
transformations (i)--(iii) also preserve the degree $\dim L_d$ of the
coefficients $E_d$ in $\Phi '$. In the ``positive'' case $l_1+...+l_r\leq n$
where $\deg q = n+1 -\sum l_a > 0$ the transformations (i)--(iii) in fact
increase degrees of the coefficients $E_d$ and are ``not allowed''. The only
exception is the operation (iii) with $f(q)=\text{const} \ q$ in the case
$l_1+...+l_r=n$ when $\deg q =1$. In Section $10$ we found the right constant
to be $-l_1!...l_r!$

\bigskip

Consider now the correlator $\Phi $,
\[ \Phi = \sum _i \frac{\Pi _a (l_a\l_i -\l'_a)}{\Pi _{j\neq i} (\l_i-\l_j)}
\ e^{\l_i (t-\t )/\hbar } \ Z_i(e^t,\hbar ) \ Z_i(e^{\t }, -\hbar ), \]
in the Calabi-Yau case $l_1+...+l_r=n+1$ (see Section $9$ for a definition
of $Z_i$).

{\bf Proposition 11.4.} {\em $(1)$ The coefficients of the power series
$Z_i(q,\hbar) = \sum _d q^d C_i(d) $ are rational functions
\[ C_i(d)=\frac{ P_d^{(i)} }{d! \hbar ^d \Pi _{j\neq i}
\Pi _{m=1}^d (\l_i -\l_j +m\hbar )} \]
where $P_d^{(i)}$ is a polynomial in $(\hbar , \l, \l')$ of degree $(n+1)d$.

$(2)$ The polynomial coefficients $E_D(p)$ in $\Phi '$
are determined by their values
\[ E_D (\l_i +d\hbar )=\Pi _a (l_a\l_i -\l'_a) \ P_d^{(i)}(\hbar )
P_{D-d}^{(i)}(-\hbar )   \]
at $p=\l_i +d\hbar $, $i=0,...,n$, $d=0,...,D$.

$(3)$ The correlators $ Z_i(q, \hbar)$
satisfy the recursion relation
\[ Z_i(q,\hbar ) = 1+ \sum_d \frac{q^d}{\hbar ^d } \frac{R_{i,d}}{d!} +
\sum _d \sum _{j\neq i}  \frac{q^d}{\hbar ^d}
\frac{ Coeff \ _i^j(d) }{\l_i-\l_j+d \hbar } \
Z_j (\frac{q}{\hbar }\frac{(\l_j-\l_i)}{d}, \frac{(\l_j-\l_i)}{d}) \]
where $R_{i,d}=R_{i,d}^{(0)}\hbar ^d + R_{i,d}^{(1)}\hbar ^{d-1} + ...$
is a polynomial of $(\hbar ,\l ,\l')$ of degree $d$, and
\[ Coeff\ _i^j(d)= \
\frac{\Pi _a \Pi _{m=1}^{l_a d} (l_a\l_i-\l'_a +m(\l_j-\l_i)/d)}
{d!\Pi _{\a \neq i}\ _{m=1}^d\ _{(\a,m)\neq (j,d)}
(\l_i-\l_{\a } +m(\l_j-\l_i)/d)}  \ .\]
For any given $\{ R_{i,d} \}$ these recursion relations have a unique solution
$\{ Z_i \} $.

$(4)$ Consider the class $\cal{P}$ of solutions $\{ Z_i \}$ to these recursion
relations which give rise to polynomial coefficients
$E_d$ in the corresponding $\Phi $. A solution from $\cal{P}$ is
uniquely determined by the first two coefficients
$R_{i,d}^{(0)}, R_{i,d}^{(1)}, i=0,...,n, d=1,...,\infty ,$
of its initial condition (that is by the first two terms in the
expansion of $Z_i= Z_i^{(0)} + Z_i^{(1)}/\hbar + ... $
as power series in $1/\hbar $).

$(5)$ The class $\cal{P}$ is invariant with respect to the following
operations:

(a) simultaneous multiplication $Z_i \mapsto f(q) Z_i $ by a power series
of $q$ with $f(0)=1$;

(b) changes $ Z_i (q,\hbar) \mapsto e^{ \l_i  f(q)/\hbar }
Z_i(q e^{f(q)},\hbar ) $ with $f(0)=0$;

(c) multiplication $Z_i \mapsto \exp (C f(q)/\hbar ) Z_i$ where $C$
is a linear function of $(\l, \l')$ and $f(0)=0$.}

\bigskip

{\em Proof.} $(3)$ We have

\[ Z_i =1+ \sum _{d>0} q^d [ \sum _{k=0}^{d-1} \hbar ^{-k-1} \int _{Y_{2,d}}
 E'_d e_0^*(\phi _i) (-c)^k ] + \sum _{d>0} q^d \hbar ^{-d} \int _{Y_{2,d}}
 E'_d e_0^*(\phi _i) \frac{(-c)^d}{\hbar + c} \]
where the integrals of the last sum have zero contributions from the type (ii)
fixed point components (Lemmas $9.7, 10.1$). Thus these integrals have a
recursive expression identical to those of Sections $9$ and $10$. The terms
of the double sum constitute the initial condition $\{ R_{i,d} \} $. The
recursion relations have the form of the decomposition of rational functions
of $\hbar $ (coefficients at powers of $Q=q/\hbar $) into the sum of simple
fractions in the case when degrees of numerators exceed degrees of
denominators. This proves existence and uniqueness of solutions to the
recursion relations.

$(1)$ follows directly from the form and topological meaning of the recursion
relations.

$(2)$ follows from the definition of $\Phi $ in terms of $Z_i$.

$(4)$ Perturbation theory: Suppose that two solutions from the class $\cal{P}$
have the same initial condition up to the order $(d-1)$ inclusively. Then
$(2)$ shows that corresponding $E_k$ for these solutions coincide for
$k<d$ and
the variation $\d E_d (p) $ vanishes at $p=\l_i + k\hbar $ for $0<k<d$.
This means that the polynomial $\d E_d$ is divisible by
$\Pi _j \Pi _{m=1}^{d-1} (p-\l_j -m\hbar)$. On the other hand $(1)$ and $(2)$
imply that the variation $\d R_{i,d}$ of the initial condition satisfies

\[ \d R_{i,d} (\hbar ) \ \Pi _a (l_a \l_i -\l'_a) \Pi _{j\neq i} \Pi _{m=1}^d
(\l_i -\l_j + m\hbar ) = \d E_d |_{p= \l_i +\hbar d }  \]
(since $R_{i,0}=1$) and thus $\d R_{i,d} $ is divisible by $\hbar ^{d-1}$.
Since $\d R_{i,d}$ is a degree $d$ polynomial, it leaves only the possibility
\[ \d R_{i,d} = \d R_{i,d}^{(0)} \hbar ^d + \d R_{i,d}^{(1)} \hbar ^{d-1} .\]
Thus if two class $\cal{P}$ solutions coincide in orders $\hbar ^{0},
\hbar ^{-1}$ then $\d R_{i,d} =0$, and thus the very solutions coincide.

$(5)$ The operations (a),(b),(c) give rise to the operations of type (i)-(iii)
for corresponding polynomials $E_d$. Thus it suffices to show that the
operations (a),(b),(c) transform a solution
$\{ Z_i \} $ of the recursion relations to another solution.

Consider in our recursion relation the coefficient
$\hbar ^{-d}q^d\ Coeff \ _i^j(d)$
responsible for the simple
fraction with the denominator $(\l_i -\l_j +d\hbar)$. The operations (a), (b),
(c)
cause respectively the following modifications in this coefficient:
\[ q^d\mapsto f(q)q^d/f(Q), \]
\[ q^d\mapsto q^d \exp \{ \frac{\l_i f(q)}{\hbar } + d f(q) -
\frac{d\l_j f(Q)}{(\l_j-\l_i)} \} , \]
\[ q^d\mapsto q^d \exp \{ C\frac{f(q)}{\hbar} +
C\frac{df(Q)}{\l_i-\l_j} \} , \]
where $Q=(\l_j-\l_i)q/d\hbar = q - (\l_i -\l_j + d \hbar ) q/ \hbar d $.
In the case of the change (b), additionally, the argument $q$ in $Z_j$
on the RHS of the recursion relation gets an  extra-factor
$\exp [ f(q)-f(Q) ]$.

The difference $Q-q$ and the exponents vanish at $\hbar = (\l_j-\l_i)/d$.
This means that
\[ \frac{q^d}{l_i -\l_j + d\hbar } \mapsto \frac{q^d}{\l_i-\l_j +d\hbar }
+ \ \text{terms without the pole} .\]
The latter terms give contributions to a new initial condition, while the
coefficient \newline $\hbar ^{-d} q^d\ Coeff \ _i^j(d)$ does not change.
It is easy to see that the required properties of the initial condition
(that the degree of $R_{i,d}(\hbar ) $ does not exceed $d$ and $R_{i,0} =1$)
are also satisfied under our assumptions about $f$ (for those contributions
involve $\hbar $ only in the combination $q/\hbar $).

\bigskip

Let us consider now the hypergeometric series

\[ Z_i^* = \sum _{d=0}^{\infty } q^d \
\frac{\Pi _{a=1}^r \Pi _{m=1}^{l_a d} (l_a\l_i -\l'_a +m\hbar)}
{\Pi _{\a =0}^n \Pi _{m=1}^d (\l_i -\l_{\a }+m\hbar)} \]
where $l_1+...+l_r=n+1$.

It is straightforward to see that $\{ Z_i^* \} $ satisfy the recursion
relations of Proposition $11.4 (3)$ (see the proof of Proposition $9.10$)
and that the formulas of Proposition $11.4 (2)$ generate corresponding

\[ \Phi ^* = \frac{1}{2\pi i} \oint \ e^{p(t-\t )/\hbar }
\sum _{d=0}^{\infty } e^{d\t }
\frac{\Pi _{a=1}^r \Pi _{m=0}^{l_a d} (l_a p - \l'_a -m\hbar)}
{\Pi _{i=0}^n \Pi _{m=0}^d (p-\l_i -m\hbar)} dp \]
with polynomial numerators. Thus $\{ Z_i \} $ is a solution from the class
$\cal{P}$.

Computation of the first two terms in the initial condition gives
\[ Z_i^*\ ^{(0)} = f(q)=\sum _{d=0}^{\infty }
\frac{(l_1d)!...(l_rd)!}{(d!)^{n+1}} \ q^d \ ,\]
\[ Z_i^*\ ^{(1)} = \l_i \sum _a l_a [g_{l_a}(q)- g_1(q)] +
(\sum _{\a } \l_{\a }) g_1(q) - \sum _a \l'_a  g_{l_a}(q) \]
where
\[ g_l=\sum _{d=1}^{\infty } q^d \frac{\Pi _a (l_a d)!}{(d!)^{n+1}}
\ (\sum _{m=1}^{ld} \frac{1}{m} ) \ .\]

\bigskip

Let us compare these initial conditions with those for $\{ Z_i \}$.

{\bf Proposition 11.5.} $Z_i^{(0)}=1, \ Z_i ^{(1)} =0$.

{\em Proof.} The first statement follows from the definition of $Z_i$ while
the second means that $\int _{Y_{2,d}} E'_d e_0^*(\phi _i) =0$ for all $d>0$.
It is due to the fact that the class $E'_d e_0^*(\phi _i)$
is a pull-back from $Y_{1,d}$.
(In fact we have just repeated an argument proving $(5)$ from Section $5$
and thus the proposition can be deduced from general properties of quantum
cohomology.)

\bigskip

Combining the last two propositions we arrive to the following

\bigskip

{\bf Theorem 11.6.} {\em The hypergeometric solution $\{ Z_i^*(q,\hbar) \} $
coincides with the solution $\{ Z_i(Q,\hbar) \} $ up to transformations
(a),(b),(c). More precisely, perform the following operations
with $\{ Z_i \}$

1) put
\[ Q=q\exp \{ \sum _a l_a [g_{l_a}(q) - g_1(q)]/f(q) \} \ ,\]

2) multiply $Z_i (Q(q),\hbar )$ by
\[ \exp \{ \frac{1}{f(q) \hbar } [\sum_{a} (l_a\l_i- \l'_a) g_{l_a} (q)
-(\sum_{\a} (\l_i-\l_{\a })) g_1(q)] \} ,\]

3) multiply all $Z_i$ simultaneously by $f(q)$.

Then the resulting functions coincide with hypergeometric
series $Z_i^*(q,\hbar)$.}

{\em Proof.} The three steps correspond to consecutive applications of
operations of type (b),(c) and (a) to $\{ Z_i \} $ and transform the
initial condition of Proposition $11.5$ to that for $\{ Z_i^*\} $.
According to Proposition $11.4$ this transforms the whole solution
$\{ Z_i \}$ to $\{ Z_i^* \} $.

\bigskip

Consider the solutions
\[ s_i = \ e^{\l_i T/\hbar } \ Z_i (e^T, \hbar ) \]
to the equivariant quantum cohomology differential equations.

{\bf Corollary 11.7.} {\em The operations

1) change $T=t+ \sum _{\a } l_a [g_{l_a} (e^t) -g_1(e^t)]/f(e^t) $,

2) multiplication by
\[ f(e^t)\exp \{ [g_1(e^t) (\sum _{\a } l_{\a }) - \sum _a \l'_a g_{l_a}(e^t)]
/(\hbar f(e^t)) \} \]
transform $\{ s_i \} $ to the hypergeometric solutions
\[ s_i^* = e^{pt/\hbar } \sum _d  e^{dt}
\frac{\Pi _a \Pi _{m=1}^{l_a d} (l_a p -\l'_a +m\hbar)}
{\Pi _{\a } \Pi _{m=1}^d (p -\l_{\a } +m\hbar)} \ |_{p=\l_i } \]
of the differential equation
\[ \Pi _{\a} (\hbar \frac{d}{dt} -\l_{\a } ) s^* =  e^t \
\Pi_a \Pi_{m=1}^{l_a} (\hbar l_a \frac{d}{dt} -\l'_a +m\hbar ) \ s^* \ .\]
For $\l'=0 ,\ \l_0+...+\l_n=0$ the solutions $s_i^*$ have the following
integral representation:
\[ \int _{\G ^n\subset \{ F_0(u)=e^t \} }
\frac{u_0^{\l_0}...u_n^{\l_n} \ du_0\w ... \w du_n }
{ \ F_1(u) \ ... \ F_r(u) \ dF_0} \]
where
\[ F_1=(1-u_1-...-u_{l_1}), \ F_2=(1-u_{l_1+1}-...-u_{l_1+l_2}), \ ...,
F_r=(1-u_{l_1+...+l_{r-1}+1}-...-u_{l_1+...+l_r}) \ \]
and $F_0=u_0...u_n$. }

\bigskip

{\bf Corollary 11.8.} {\em The hypergeometric class $S^*(t,\hbar )\in
H^*(\CC P^n)=\CC [P]/(P^{n+1})$,
\[ S^*= e^{Pt/\hbar } \sum_d e^{dt} \frac{\Pi_a \Pi_{m=0}^{l_ad}
(l_aP+m\hbar )}{\Pi_{m=1}^d (P+m\hbar)^{n+1} } \]
whose $n+1-r$ non-zero components are solutions to the Picard-Fuchs equation
\[ (\frac{d}{dt})^{n+1-r} s^* = l_1...l_r e^t \Pi_a \Pi_{m=1}^{l_a-1}
(l_a \frac{d}{dt} + m) s^* \]
for the integrals
\[ \int _{\g ^{n-r} \subset X_t'} \frac{du_0\w ...\w du_n}
{dF_0\w dF_1\w ... \w dF_r} \ ,\]
(here $X_t'=\{ (u_0,...,u_n) | F_0(u)=e^t, F_1(u)=0, ..., F_r(u)=0 \} $)
are obtained from the class $S$ (describing the quantum cohomology
$\cal D$-module for the Calabi-Yau complete intersection
$X^{n-r}\subset \CC P^n$),
\[ S=e^{PT/\hbar } \sum_d e^{dT} (e_0)_* (\frac{E_d}{\hbar+c_1^{(0)}}) ,\]
by the change
\[ T=t +  \sum _a l_a [g_{l_a}(e^t)-g_1(e^t)]/f(e^t) \]
followed by the multiplication by $f(e^t)$.}

{\em Proof.} Corollary $11.7$ shows that for $\l'=0, \sum \l_{\a }=0$
these change and multiplication transform the corresponding equivariant
classes $S'$ and $S'\ ^*$ to one another. The class $-p$ in the formula for
$s_i^*$ in Corollary $11.7$ is the equivariant Chern class of the Hopf line
bundle over $\CC P^n$. In the limit $\l =0$ it becomes $-P$ while $S'$ and
$S'\ ^*$ transform to their non-equivariant counterparts $S$ and $S^*$.

{\em Remarks.} 1) Notice that the components $S_0^*$ and $S_1^*$ in
\[ S^*=l_1...l_r[P^r S_0^* (t)+ P^{r+1} S_1^*(t) + ... + P^n S_n^*(t)] \]
are exactly $f(e^t)$ and $tf(e^t)+\sum_a l_a [g_{l_a }(e^t) - g_1(e^t)]$
respectively. Thus the inverse transformation from $S^*$ to $S$ consists in
division by $S_0^*$ followed by the change $T= S_1^*(t)/S_0^*(t)$ in complete
accordance with the recipe \cite{COGP, BVS, HG1} based on the mirror
conjecture.

2) According to \cite{B} the $(n-r)$-dimensional manifolds $X_t'$ admit a
Calabi-Yau compactification to the family $\bar{X}_t'$ of
{\em mirror manifolds } of the Calabi-Yau complete intersection
$X^{n-r}\subset \CC P^n$. The Picard-Fuchs differential equation from
Corollary $11.8$ describes variations of complex structures for $\bar{X}'$.
This proves the mirror conjecture (described in detail in \cite{BVS})
for projective Calabi-Yau complete intersections and confirms the enumerative
predictions about rational curves and quantum cohomology algebras made there
(and in some other papers) on the basis of the mirror conjecture.

3) The description \cite{AM} of the
quantum cohomology algebra of a Calabi-Yau $3$-fold in terms of the numbers
$n_d$ of rational curves of all degrees $d$ (see for instance \cite{HG1} for
the description of the corresponding class $S$ in these terms) has been
rigorously justified in \cite{M}. Combining these results with Corollary
$11.8$ we arrive to the theorem formulated in the introduction.

\bigskip

{\em Proof of The Main Lemma.}

In our construction of the map $\m: L_d\to L'_d$ we will denote $L_d$
the moduli space of stable maps $C\to \CC P^n\times \CC P^1$ of bidegree
$(d,1)$ with no marked points (it also has dimension $d(n+1)+n$).
The construction works for any given number
of marked points but produces a map which is the composition of $\m $ with
the forgetful map. In this form it applies to the submanifold of stable maps
with two marked points confined over $0$ and $\infty $ in $\CC P^1$ (this
submanifold is what we denoted $L_d$ in the formulation of The Main Lemma).

Let $\psi: C\to \CC P^n \times \CC P^1 $ be a stable genus $0$ map of bidegree
$(d,1)$. Then $C=C_0 \cup C_1 ... \cup C_r$ where $C_0$ is isomorphic to
$\CC P^1$
and $\psi | C_0$ maps $C_0$ onto the graph of a degree $d'\leq d$ map $\CC P^1
\to \CC P^n$, and for $i=1,...,r$ the bidegree $(d_i,0)$ map $\psi | C_i$
sends $C_i$ into the slice $\CC P^n \times \{ p_i \} $ where $p_i \neq p_j$
and $d_1+...+d_r=d-d'$.

The map $\m :L_d \to L'_d$ assigns to $[\psi ]$ the $(n+1)$-tuples
$(f_0 g : f_1 g : ... : f_n g)$ of polynomials ($=$ binary forms) on
$\CC P^1$ where
$g$ is the polynomial of degree $d-d'$ with roots $(p_1,...,p_r)$ of
multiplicities $(d_1,...,d_r)$ and the tuples $(f_0:...:f_n)$ of degree $d'$
polynomials (with no common roots, including $\infty $)
is the one that describes the map $\psi | C_0$.

In order to prove that the map $\m $ is regular let us give it another, more
invariant description.

Denote $\hat{L}_d$ the moduli space of bidegree $(d,1)$ stable maps with
an extra-marked point and pull back to $\hat{L}_d$ the line bundle
\[ H:= Hom (\pi _1^* \calo _{\CC P^n} (1), \pi _2^* \calo _{\CC P^1} (d)) \]
by the evaluation map $e: \hat{L}_d \to \CC P^n\times \CC P^1$ (where $\pi _i$
are projections to the factors). Consider the push-forward sheaf
$H^0:=R^0\pi _* e^* (H)$
of the locally free sheaf $e^* H$ along the forgetful map $\pi :\hat{L}_d\to
L_d$. To a small neighborhood $U\subset L_d$, it assigns the $\calo _U$
-module $H^0(\pi ^{-1}(U), e^* H)$ of sections of $e^* H$.

{\em Claim.} {\em $1$) $H^0$ is a rank $1$ locally free sheaf on $L_d$.

$2$) The fiber at $[\psi ]$ of the corresponding line bundle can be
identified with
\[ H^0(C_0, (\psi | C_0)^*(H) \otimes \calo (-[p_1])^{\otimes d_1} ...
\otimes \calo (-[p_r])^{\otimes d_r}) .\]

$3$) The kernel of the natural map
\[ h: H^0(C, \psi ^*\pi _1^*(\calo_{\CC P^n}(1))) \to
   H^0(C, \psi ^*\pi _2^*(\calo_{\CC P^1}(d))) =
   H^0(\CC P^1, \calo (d)) \]
defined by a nonzero vector in this fiber consists of the sections
vanishing identically on $C_0$.}

 Using this, we pick $n+1 $ independent sections of $\calo_{\CC P^n}(1)$
(that is homogeneous coordinates on $\CC P^n$), define corresponding
sections of $e^* \pi _1^* \calo_{\CC P^n} (1)$ and apply the map $h$.
By this we obtain a degree $1$ map from the total space of the line bundle
$H^0$ to the linear space $\CC ^{n+1} \otimes H^0(\CC P^1, \calo (d))$. Since
the homogeneous coordinates on $\CC P^n$ nowhere vanish simultaneously,
we obtain a natural map
\[  L_d \to L'_d = Proj (\CC ^{n+1}\otimes H^0(\CC P^1, \calo (d))) \]
which sends $[\psi ]$ to $(f_0 g:...:f_n g)$ and conclude that $\m $ is
regular.

The remaining statements of The Main Lemma are proved by looking at
localizations of the equivariant class $p$ at the $S^1\times T^{n+1}$-fixed
points in $L'_d$ and $L_d$ (in this paragraph we use the notation $L_d$ for
the same space as in the formulation of The Main Lemma).
The fixed points in $L'_d$ are represented by
the vector-monomials $(0:...:0:x^{d'}:0:...:0)$ where $p$ localizes to
$\l_i + d'\hbar $. A fixed point in $L_d$ is represented by $\psi $ with
$\psi (C_0) = (0: ... :0:1: 0:...:0)$, $r=2$, $p_0=0$, $p_1=\infty $ and
the maps $\psi | C_k : C_k \to \CC P^n$, $k=1,2$ representing $T^{n+1}$-fixed
points respectively in $Y_{2,d'}$ and $Y_{2,d-d'}$ such that their (say)
second marked points are mapped to the point $\psi (C_0)$.
This implies that the class $\m^*(p)$ localizes to $\l_i+d'\hbar $ at such a
fixed point and thus the pull back of $p$ to the fixed point set
\[ \{ [\psi ]\in Y_{2,d'}\times Y_{2,d-d'} | (e_2\times e_2) ([\psi ])
\in \D \subset Y\times Y \]
of the $S^1$-action on $L_d$ coincides with the pull back through the common
marked point of the $T^{n+1}$-equivariant class $p+d'\hbar$ on the diagonal
$\D = \CC P^n$. Now localizations of $\int _{L_d} e^{p(t-\t )} {\cal E} _d $
to the fixed points of $S^1$-action identify the form of the correlator
$\Phi $ given in The Main Lemma with the definition of $\Phi $ as
the convolution of $S'(t,\hbar )$ and $S'(\t ,-\hbar )$.

\bigskip

In order to justify the {\em claim} we need to compute the space of global
sections of the sheaf $e^* (H)$ over the formal neighborhood of the fiber
$\pi ^{-1} ([\psi ])$ of the forgetful map $\pi : \hat{L}_d \to L_d$.
The fiber itself is isomorphic to the tree-like genus $0$ curve $C$.
Let $(x_j,y_j), j=1,...,N\geq r$ be some local parameters on irreducible
components of $C$ near the singular points such that $\e_j = x_jy_j$ are local
coordinates on the {\em orbifold} $L_d$ near $[\psi ]$ (one should add
some local coordinates $\e'$ on the stratum $\e_1=...\e_N=0$ of stable maps
$C\to \CC P^n$ in order to construct a complete
local coordinate system on $L_d$). Such a description of local coordinates on
$L_d$ follows from the very construction of the moduli spaces of stable maps
to ample manifolds; we refer the reader to \cite{Kn, BM} for details.

A line bundle over the neighborhood of $C\subset \hat{L}_d$ can be specified
by the set
\[ u_j (x_j^{\pm 1}, \e ),
v_j (y_j^{\pm 1}, \e ), \ j=1,...,N, \]
of non-vanishing functions describing transition maps between trivializations
of the bundle inside and outside the neighborhoods
(with local coordinates $(x_j,y_j,\e_1,...,\hat{\e _j},...,\e_N,\e')$)
of the double points.

Let us consider first the following model case. Suppose that $C$
consists of $r+1$ irreducible components $(C_0, C_1,...,C_r)$ such that each
$C_j$ with $j>0$ intersects $C_0$ at some point $p_j$. Let $x_j$ be
the local parameter on $C_0$ near $p_j$, and the line bundle
(of the degree $-d_j\leq 0$ on $C_j$) be specified by $v_j=y_j^{-d_j}$.

In the neighborhood of $p_j$ a section of such a bundle is given
by a function $s(x_j,y_j,\hat{\e _j})$ satisfying
\[ s=y_j^{-d_j}s_j(y_j^{-1}, \e) \]
where the function $s_j$ represents the section in the trivialization
over the neighborhood of $C_j-p_j$. Here $\hat{\e_j}$ means that $\e_j$
is excluded from the set of coordinates $\e $ (remember that $\e_j=x_jy_j$).
This implies that $s_j = \e _j^{d_j} f_j(y_j^{-1}\e_j, \e)$ where $f_j$
is some regular function. Thus this section in the neighborhood of $p\in C_0$
is given by a function $ s(x_j, \e ) = x_j^{d_j}f_j(x_j,\e)$ with zero of
order $d_j$ at $x_j=0$, and the restriction of this section to the
neighborhood of $C_j$ is determined by $s$.

In other words, the $\CC [[\e ]]$-module of global sections in the
formal neighborhood of $C$ identifies with the module of global sections
on $C_0$ for the line bundle given by the loops $x_j^{-d_j} u_j$
instead of $u_j$ (this corresponds to the subtraction of the divisor
$\sum d_j [p_j]$.

The more general situation where $v_j$ is the product of $y_j^{-d_j}$ with
an invertible function $w_j(y_j, x_j ,\hat{\e_j })$
preserves the above conclusion with
$w^{-1}s=x_j^{d_j} f_j(x_j,\e )$ instead of $s$.

Obviously, the above computation bears dependence on additional parameters.

Now we apply our model computation to the neighborhood of a general tree-like
curve $C$ {\em inductively} by decomposing the tree into simpler ones
starting from
the root component $C_0$. We conclude that the $\CC [[\e]] $-module of
sections of the bundle $e^*(H)$ is identified with the module of sections
of some line bundle
over the product of $C_0$ with the polydisk with coordinates
$(\e_1,...,\e_r, ..., \e_N, \e')$,
and that this line bundle is $e^*(H)$ for $C_0$
(given by the loops $u_j$ in our current
notations) twisted by the loops $x_j^{-d_j}$
in the punctured neighborhoods of the points $(p_1,...,p_r)$, where
$(d_1,...,d_r)$ are the degrees of the maps $\psi |C_j: C_j \to \CC P^n$
(in the notations of the {\em claim} so that $d_1+...+d_r=d-d'$).

This implies that the $\CC [[\e ]]$-module ${\cal H}^0$
of global sections can be identified with the
module of those global sections of the degree $d-d'$ locally free sheaf
$(\psi |C_0)^* (H) \otimes \CC [[\e ]]$ which have zeroes of order
$d_j$ at $p_j$ for $j=1,...,r$. In particular

1) ${\cal H}^0$ is a free $\CC [[\e ]]$-module of rank $1$,

2) ${\cal H}^0 \otimes _{\CC [[\e ]] } (\CC [[\e ]]/(\e))$ is the
$1$-dimensional space $H^0|_{[\psi]} $ described in the {\em claim}, and

3) non-zero vectors in $H^0|_[\psi]$ represent sections of $\psi ^*(H)$ over
$C$ non-zero on $C_0$ (and thus their product with a non-zero on $C_0$
section of
$\psi ^*\pi _1^*(\calo _{\CC P^n} (1))$ can not vanish identically on $C_0$.)

Factorization by the discrete group  $Aut (\psi )$ preserves $(1-3)$ with
$\CC [[\e ]]$ replaced by $\CC [[\e ]]^{Aut (\psi )}$.

\end{document}